\documentclass{amsart}
\usepackage{amsmath,amsfonts,amsthm,amscd}
\newtheorem{thm}{Theorem}[section]
\newtheorem{pro}[thm]{Proposition}
\newtheorem{cor}[thm]{Corollary}
\newtheorem{lem}[thm]{Lemma}
\newtheorem{claim}[thm]{Claim}

\newtheorem{defn}[thm]{Definition}
\newtheorem{quot}[]{Theorem}
\newcommand{\dime}{\operatorname{dim}}
\newcommand{\Obj}{\operatorname{Obj}}
\newcommand{\Mor}{\operatorname{Mor}}
\newcommand{\chan}{\operatorname{char}}
\newcommand{\degr}{\operatorname{deg}}
\newcommand{\expo}{\operatorname{exp}}
\newcommand{\spant}{\operatorname{span}}
\newcommand{\bone}{\hat{\beta}_1}
\newcommand{\btwo}{\hat{\beta}_2}

\begin{document}
\title{Cohomology of uniformly powerful $p$-groups}

\author{William Browder}
\author{Jonathan Pakianathan}

\begin{abstract}
In this paper we will study the cohomology of a
family of $p$-groups associated to $\mathbb{F}_p$-Lie algebras.
More precisely we study a category $\mathbf{BGrp}$ of $p$-groups
which will be equivalent to the category of $\mathbb{F}_p$-bracket
algebras (Lie algebras minus the Jacobi identity). We then show that
for a group $G$ in this category, its $\mathbb{F}_p$-cohomology
is that of an elementary abelian $p$-group if and only if it is
associated to a Lie algebra.

We then proceed to study the exponent of $H^*(G ;\mathbb{Z})$ in
the case that $G$ is associated to a Lie algebra $\mathfrak{L}$. To do this,
we use the Bockstein spectral sequence and derive a formula that
gives $B_2^*$ in terms of the Lie algebra cohomologies of
$\mathfrak{L}$. We then expand some of these
results to a wider category of $p$-groups. In particular we calculate
the cohomology of the $p$-groups $\Gamma_{n,k}$ which are defined to
be the kernel of the mod $p$ reduction
$
GL_n(\mathbb{Z}/p^{k+1}\mathbb{Z}) \overset{mod}{\longrightarrow}
GL_n(\mathbb{F}_p).
$

\end{abstract}

\maketitle

\section{Introduction and Motivation}

Throughout this paper, $p$ will be an odd prime. First some definitions:

\begin{defn}
\label{defn: omega1}
Given a $p$-group $G$, $\Omega_1(G) = \langle g \in G : g^p=1\rangle$
where the brackets mean ``smallest subgroup generated by''.
$G$ is called $p$-central if $\Omega_1(G)$ is central.
\end{defn}
\begin{defn}
\label{defn: ppower}
For any positive integer $k$, $G^{p^k} = \langle g^{p^k} : g \in G \rangle$.
\end{defn}
\begin{defn} 
$Frat(G) = G^p[G,G].$
\end{defn} 

We will study a category $\mathbf{BGrp}$ 
of $p$-groups which is naturally equivalent to the category of 
bracket algebras over $\mathbb{F}_p$. This is exactly the category of 
$p$-central, $p$-groups $G$ which have $\Omega_1(G)=G^p=Frat(G)$. 
For such a group $G$,
the associated bracket algebra will be called $Log(G)$.

In the case that one of these groups is associated to a Lie algebra, we
will show that it has the same
$\mathbb{F}_p$-cohomology as an elementary abelian $p$-group. More precisely we
will show: (The necessary definitions for the theorems quoted in this introductory section can be found in the relevant parts of the paper.)

\begin{quot}[~\ref{thm: Fp cohomology}]
Let $G \in \Obj(\mathbf{BGrp})$ and $n=\dime(\Omega_1(G))$. Then 
$$
H^*(G;\mathbb{F}_p) = \wedge(x_1,\dots,x_n) \otimes 
\mathbb{F}_p[s_1,\dots,s_n]
$$
(where the $x_i$ have degree 1 and the $s_i$ have degree 2) if and only
if $Log(G)$ is a Lie algebra. When this is the case, the polynomial
algebra part restricts isomorphically to that of $H^*(\Omega_1(G);\mathbb{F}_p)$
and the exterior algebra part is induced isomorphically from that of
$H^*(G/\Omega_1(G);\mathbb{F}_p)$ via the projection homomorphism.
\end{quot}
\noindent 
(This theorem was proven independently by T. Weigel for the case $p \neq 3$
in \cite{W1}, \cite{W2}.)

However the integral cohomology of these groups is complicated and indeed
we will show that one can recover the Lie algebra associated to the
group from knowledge of the Bockstein on its $\mathbb{F}_p$-cohomology.
To do this we calculate the full structure of the $\mathbb{F}_p$-cohomology
as a Steenrod-module. As an application, in corollary~\ref{cor: comodule},
we completely determine the comodule algebra structure (see \cite{W2}) of 
$H^*(G; \mathbb{F}_p)$.

With the mild additional hypothesis that the associated Lie algebra lifts to
one over $\mathbb{Z}/p^2\mathbb{Z}$, we compute $B_2^*$ of the Bockstein 
spectral sequence in terms of the Lie algebra cohomologies of the 
corresponding Lie algebra:

\begin{quot}[~\ref{thm: B2}]
Let $G \in \Obj(\mathbf{BGrp})$ with $Log(G)=\mathfrak{L}$, a Lie algebra, 
and suppose that $\mathfrak{L}$ lifts to a Lie algebra over 
$\mathbb{Z}/p^2\mathbb{Z}$.
Then $B_2^*$ of the Bockstein spectral sequence for $G$ is given by
$$
B_2^*=\oplus_{k=0}^{\infty}H^{*-2k}(\mathfrak{L};S^k).
$$
\end{quot} 
Here $S^k$ is the $\mathfrak{L}$-module of symmetric $k$-forms on $\mathfrak{L}$ with the usual action. (This action is described in detail before the proof of the theorem.) 
	
The integral cohomology of $p$-groups is a rich subject and
there are many theorems and conjectures about the exponent of $\Bar{H}^*(P ;
\mathbb{Z})$. (Where the bar denotes reduced cohomology and we recall
that the exponent of an abelian group $G$ is the smallest positive integer $n$ such that
$ng=0$ for all $g \in G$). We will obtain some partial results on the exponent 
of the integral cohomology of the $p$-groups
studied in this paper via the Bockstein spectral sequence.

In the last section of the paper, we extend the results mentioned to a
bigger family of groups: the uniform, $p$-central, $p$-groups.
(Elementary abelian $p$-groups are uniform and more generally one defines
inductively that a $p$-central, $p$-group $G$ is uniform if and only if 
$\Omega_1(G)=G^{p^k}$ for some nonnegative integer $k$ and $G/\Omega_1(G)$ is 
itself uniform. Thus a uniform $p$-group will give rise to a tower of 
uniform $p$-groups
called a uniform tower where each group in the tower is the quotient of the previous group $G$ by $\Omega_1(G)$. In this paper such uniform towers will be
indexed so that $G_1$ is always an elementary abelian $p$-group. When this
is done, $G_2$ will always correspond to a group in the category 
$\mathbf{BGrp}$.)

Specifically we prove:

\begin{quot}[~\ref{thm: main}]
Fix $p \geq 5$. Let
$$
G_{k+1} \rightarrow G_k \rightarrow \dots \rightarrow G_2 \rightarrow
G_1 \rightarrow 1
$$
be a uniform tower with $k \geq 2$. Let $\mathfrak{L}=Log(G_2)$
and let $c_{ij}^k$ be the structure constants of $\mathfrak{L}$ with
respect to some basis. Then for suitable
choices of degree 1 elements $x_1,\dots,x_n$ and degree 2
elements $s_1,\dots,s_n$, one has
$$
H^*(G_k ; \mathbb{F}_p) \cong \wedge^*(x_1,\dots,x_n) \otimes
\mathbb{F}_p[s_1,\dots,s_n]
$$
with 
\begin{align*}
\begin{split}
\beta(x_t) =& -\sum_{i < j}^n c_{ij}^tx_ix_j \\
\beta(s_t) =& \sum_{i,j=1}^n c_{ij}^ts_ix_j + \eta_t
\end{split}
\end{align*}
for all $t = 1,\dots n$. Furthermore the $\eta_t$ define a 
cohomology class $[\mathbf{\eta}] \in H^3(\mathfrak{L} ; ad)$ which vanishes if
and only if there exists a uniform tower
$$
G_{k+2}' \rightarrow G_{k+1}' \rightarrow G_k \rightarrow
\dots \rightarrow G_2 \rightarrow G_1 \rightarrow 1
$$ or in other words if and only if $Log(G_k)$ (which is a Lie algebra over
$\mathbb{Z}/p^{k-1}\mathbb{Z}$), has a lift to a Lie algebra
over $\mathbb{Z}/p^k\mathbb{Z}$. Thus in this case where the Lie algebra lifts,
one can drop the $\eta_t$ term in the formula for the Bockstein.
\end{quot}

A version of this theorem for the category of powerful, $p$-central, $p$-groups
is also stated in theorem~\ref{thm: last}.
 
Some examples of $p$-groups covered in this paper are the groups
$\Gamma_{n,k}$ for $n,k \geq 1$, where the group $\Gamma_{n,k}$ is defined
as the kernel of the map
$$
GL_n(\mathbb{Z}/p^{k+1}\mathbb{Z}) \overset{mod}{\longrightarrow}
GL_n(\mathbb{F}_p).
$$

\section{Bracket groups}
\subsection{Some preliminaries}
\label{sec: prelim}

First we will recall a method of obtaining a bracket algebra
from a $p$-group under suitable conditions.

Let us consider the following situation. Suppose we have a central short
exact sequence of finite groups:
$$
1 \rightarrow V \rightarrow G \overset{\pi}{\rightarrow} W \rightarrow 1
$$
where W and V are elementary abelian $p$-groups, then we define 
\begin{align}
 \langle \cdot,\cdot\rangle  : W \times W  \rightarrow V 
 \text{ given by }  
 \langle  x,y\rangle =\Hat{x}\Hat{y}\Hat{x}^{-1}\Hat{y}^{-1}
\end{align} 
where $\Hat{x},\Hat{y}$
are lifts of $x,y$ to $G$, i.e.,\ $\pi(\Hat{x})=x,\pi(\Hat{y})=y$.

We also define a $p$-power map function 
\begin{align}
\phi : W \rightarrow V \text{ given by }
\phi(x) = \Hat{x}^{p}.
\end{align}
It is routine to verify that $\langle \cdot, \cdot\rangle$ is a well-defined
alternating, bilinear map and that $\phi$ is a well-defined linear map. 
(the linearity of $\phi$ uses that $p$ is odd.)

We will now restrict ourselves to the case where the $p$-power map $\phi$
is an isomorphism. So $V$ and $W$ are isomorphic.
In this case one can define 
\begin{align}
[\cdot,\cdot] : W \otimes W  \rightarrow W \text{ by } 
 [w_1,w_2] = \frac{1}{2}\phi^{-1}(\langle w_1,w_2\rangle ).
\end{align} 
Then it is easy to see this
is still an alternating, bilinear map on $W$. Note $W$ is a vector space
over $\mathbb{F}_p$, the field on $p$ elements. (The factor of $\frac{1}{2}$
is put in for convenience in order to avoid messy expressions later on.)

\begin{defn}
A bracket algebra over $\mathbb{F}_p$ is a finite dimensional vector space
 $W$ over $\mathbb{F}_p$ equipped with an alternating, bilinear form
$[\cdot,\cdot] : W \otimes W \rightarrow W$. The dimension of a bracket algebra is the dimension of the underlying vector space.
\end{defn}
\begin{defn}
$\mathbf{Brak}$ is the category
whose objects consist of bracket algebras over $\mathbb{F}_p$ and with morphisms the linear maps which preserve the brackets.
\end{defn}

Now we can define a category of $p$-groups which is essentially equivalent
to the category $\mathbf{Brak}$. These groups will
be our primary objects of study in the first part of this paper.

\begin{defn}
A $p$-power exact extension is a central short exact sequence:
$$
1 \rightarrow V \rightarrow G \overset{\pi}{\rightarrow} W \rightarrow 1
$$
where V, W are elementary abelian and where the
$p$-power map $\phi$ is an isomorphism. 
\end{defn}
Note in such an extension, $V=\Omega_1(G)$ and $W=G/\Omega_1(G)$ (see definition~\ref{defn: omega1}).
So any homomorphism between two groups $G_1,G_2$ in the middle of 
$p$-power exact extensions takes $V_1=\Omega_1(G_1)$
into $V_2=\Omega_1(G_2)$ and hence induces maps $V_1 \rightarrow V_2$ and
$W_1 \rightarrow W_2$ where $W_i=G_i/\Omega_1(G_i)$. These will be called respectively the induced maps on the $V, (W)$ level. 
\begin{defn}
$\mathbf{BGrp}$ is the category whose objects consist of
finite groups $G$ which fit in the middle of a $p$-power exact extension.
These will be called bracket groups. The morphisms in $\mathbf{BGrp}$ are just
the usual group morphisms between these groups, 
except we will identify two morphisms if
they induce the same map on the $V$ and $W$ levels. 
\end{defn} 

One can show that the two categories $\mathbf{Brak}$ and $\mathbf{BGrp}$ are naturally equivalent (see \cite{pak}). This follows from a general exponent-log
correspondence (see \cite{W1}, \cite{W3} or \cite{laz}), suitably reworded for
our purposes. We will give a short
description of the functors involved in this equivalence of categories.

Define a covariant functor
$Log : \mathbf{BGrp} \rightarrow \mathbf{Brak}$ as follows, to a bracket group
$G$ we associate the bracket algebra $Log(G)=(W, [\cdot,\cdot])$
which is obtained as explained before. Notice that the underlying vector
space of $Log(G)$ is just $G/\Omega_1(G)$ so given $\psi \in \Mor(G_1,G_2)$ we
note $\psi$ induces a well-defined linear map 
$Log(\psi) : Log(G_1) \rightarrow Log(G_2)$,
and it is easy to check that this is a map of bracket algebras.

A description of the inverse functor $Exp : \mathbf{Brak} \rightarrow
\mathbf{BGrp}$ is given as follows.  
Given $(L, [\cdot,\cdot])$ a bracket algebra, let $K$ be a free $\mathbb{Z}/
p^2\mathbb{Z}$-module of rank equal to the dimension of $L$. Let $\phi :
L \rightarrow K$ and $\pi : K \rightarrow L$ be injective/surjective maps
such that $\phi \circ \pi = p$ (multiplication by p). Then $Exp(L)=(K, \circ)$
where 
\begin{align}
\begin{split}
\label{eq: Expdef}
l \circ m = l + m + \phi[\pi(l),\pi(m)].
\end{split}
\end{align} 
(Here $+$ is the addition of
$K$ as a $\mathbb{Z}/p^2\mathbb{Z}$-module. Notice in this notation, $0$ is the
identity of $Exp(L)$, $l^{-1}$ corresponds to $-l$, $l^p$ corresponds to
$pl$, and the $p$-power map is indeed $\phi$. )

For clarity we state the following proposition which summarizes these facts:

\begin{pro}
The functors $Log : \mathbf{BGrp} \rightarrow \mathbf{Brak}$ and 
$Exp : \mathbf{Brak} \rightarrow \mathbf{BGrp}$ give a natural equivalence
between the categories $\mathbf{BGrp}$ and $\mathbf{Brak}$. Thus to every
bracket algebra there naturally corresponds a unique bracket group.
\end{pro}
	
\noindent
\textbf{Example:}

A direct computation shows that the group $\Gamma_{n,2}$ which is the kernel of the map
$$
GL_n(\mathbb{Z}/p^3\mathbb{Z}) \overset{mod}{\longrightarrow}
GL_n(\mathbb{F}_p),
$$
is a bracket group, and that $Log(\Gamma_{n,2}) = \mathfrak{gl}_n$,
the Lie algebra of $n \times n$ matrices.

Now we will study the groups in $\mathbf{BGrp}$ from
a cohomological viewpoint. 
Recall if $G \in \Obj(\mathbf{BGrp})$ then $G$ fits in a central short
exact sequence
$$
1 \rightarrow V \rightarrow G \overset{\pi}{\rightarrow} W \rightarrow 1
$$
where we can identify $V$ and $W$ via the p-power isomorphism $\phi$. 
If we do not
put the restriction that $\phi$ is an isomorphism then such extensions as
above are in one to one correspondence with $H^2(W;V)$. Now for
$W$ an elementary abelian p-group of rank $n$ (recall $p$ is odd), we have that
$$
H^*(W;\mathbb{F}_p) \cong \wedge^*(x_1,\dots,x_n) \otimes \mathbb{F}_p[
\beta{x_1},\dots,\beta{x_n}]
$$
the tensor product of an exterior algebra on degree 1 generators and
a polynomial algebra on degree 2 generators which are the Bocksteins
of the degree 1 generators. (Here $\beta: H^1(W; \mathbb{F}_p) \rightarrow
H^2(W;\mathbb{F}_p)$ is the Bockstein. )
In basis free notation, this can be written
$$
H^*(W;\mathbb{F}_p) \cong \wedge^*(W^*) \otimes \mathbb{F}_p[\beta(W^*)]
$$
where $W^*=H^1(W;\mathbb{F}_p)$ is the dual vector space of $W$ and 
$\beta(W^*)$ is its image under the Bockstein $\beta
.$
Now,
\begin{align}
\begin{split}
\label{eq: silly}
H^2(W;V) \cong & \ H^2(W;\mathbb{F}_p) \otimes V  \text{ by the Universal
Coefficient Theorem} \\
 \cong & \ (H^2(W;\mathbb{F}_p))^n  \text{  using a choice of basis of } V.
\end{split}
\end{align}
Note that a bracket on $W$ gives us a map $br: \wedge^2(W) \rightarrow W$.
If we take the dual of this map we get a map
$br^*: W^* \rightarrow \wedge^2(W^*)$. (Once we identify the dual of 
$\wedge^2(W)$ with $\wedge^2(W^*)$ in the usual way.)
With this notation, it is
easy to argue by comparisons (see \cite{BeC}) and from the description of
the $Exp$ functor in equation~\ref{eq: Expdef} that the extension 
elements are of the form
$ \beta(x_i) + br^*(x_i) $ for $i=1,\dots,n$. (If we use 
a suitable basis for $V$ in
equation~\ref{eq: silly}. One chooses a basis $E$ for $W$, the
basis $\phi(E)$ for $V$ and the canonical dual basis for $V^*$
and $W^*$ so for example $\{ x_1,\dots, x_n \}$ are dual to $E$.) 
In this formulation
$br^*(x_i)$ are just the components of the bracket
of $Log(G)$, with respect to the basis $E$.
We will use this form of the extension element from now on.

\subsection{Lie Algebras and co-Lie algebras}
\begin{defn}
Given a bracket algebra $(W,[\cdot,\cdot])$ we define the Jacobi form
of $(W,[\cdot,\cdot])$ to be 
$$
J(x,y,z) = ([[x,y],z] + [[y,z],x] + [[z,x],y]).
$$
This is easily checked to be an alternating 3-form. Thus
$J : \wedge^3(W) \rightarrow W$. Bracket algebras with $J=0$ are called
Lie algebras.
\end{defn}

Given a bracket algebra, $(W, [\cdot,\cdot])$, the bracket defines a map
$br : \wedge^2(W) \rightarrow W$. Thus by taking duals and identifying
the dual of $\wedge^*(W)$ with $\wedge^*(W^*)$ in the usual way, we obtain a map
$br^* : W^* \rightarrow \wedge^2(W^*)$, where $W^*$ is the dual vector space
of $W$. This motivates the following definition:

\begin{defn}
A co-bracket algebra is a vector space $W^*$ equipped with a linear map
$br^* : W^* \rightarrow \wedge^2(W^*)$. We will always extend  
$br^*$ as a degree 1 map on the whole of the graded algebra $\wedge^*(W^*)$
in the unique way that makes it a derivation on that algebra. We call such a 
co-bracket algebra a co-Lie algebra if $br^* \circ br^* = 0$.
\end{defn}

As was mentioned before, 
the dual of a bracket algebra, has a natural
co-bracket algebra structure and similarly, the dual of a co-bracket algebra
has a natural bracket algebra structure.

Given a co-bracket algebra $W^*$, it is easy to check, 
as $br^*$ is a derivation, that $br^* \circ br^* = 0$ if 
and only if $br^* \circ br^* : W^*=\wedge^1(W^*) \rightarrow \wedge^3(W^*)$ 
is zero.
The dual of this map is a map from $\wedge^3(W) \rightarrow W$ (where we
are calling $W^{**}=W$). It is easy to check that this is just the Jacobi
form defined before for the bracket algebra $W$. Thus we can conclude that
the dual of a co-Lie algebra is a Lie algebra and vice versa.

Let $\pi : \wedge^2(W) \rightarrow W \otimes W$ be the map
defined by $\pi(a \wedge b) = a \otimes b - b \otimes a$. 
Thus the usual definition of a module over a bracket (Lie) algebra can be 
worded in the following
more categorical way:

\begin{defn}
We say $V$ is a module over the bracket algebra $(W, br)$ if we
have a map $\lambda : W \otimes V \rightarrow V$ such that the following
diagram commutes:
$$
\begin{CD}
W \wedge W \otimes V @>> br \otimes 1 > W \otimes V   \\
@V (1 \otimes \lambda) \circ (\pi \otimes 1) VV @VV \lambda V\\
W \otimes V @> \lambda >> V 
\end{CD}
$$
\end{defn}

Taking the dual of this definition, and noting that the dual of $\pi$,
$\pi^* : W^* \otimes W^* \rightarrow \wedge^2(W^*)$ is just the 
canonical quotient map,
 gives us the following definition:

\begin{defn}
We say $V^*$ is a comodule over the co-bracket algebra $(W^*, br^*)$ if we
have a map $\lambda^* : V^* \rightarrow W^* \otimes V^*$ such that the following
diagram commutes:
$$
\begin{CD}
W^* \wedge W^* \otimes V^* @<< br^* \otimes 1 < W^* \otimes V^*\\
@A 1\wedge \lambda^* AA @AA \lambda^* A\\
W^* \otimes V^* @< \lambda^* << V^*
\end{CD}
$$
\end{defn}

It is easy to see that the dual of a comodule over the co-bracket algebra
$W^*$ is a module over the bracket algebra $W=W^{**}$ and vice-versa.

\subsection{The $\mathbb{F}_p$-cohomology of a group in $\mathbf{BGrp}$}
In this section we will show that a group in $\mathbf{BGrp}$ has the
$\mathbb{F}_p$-cohomology of an elementary abelian $p$-group if and only
if it is associated to a Lie algebra. More precisely we will prove:

\begin{thm}
\label{thm: Fp cohomology}
Let $G \in \Obj(\mathbf{BGrp})$ and $n=\dime(\Omega_1(G))$. Then 
$$
H^*(G;\mathbb{F}_p) = \wedge(x_1,\dots,x_n) \otimes 
\mathbb{F}_p[s_1,\dots,s_n]
$$
(where the $x_i$ have degree 1 and the $s_i$ have degree 2) if and only
if $Log(G)$ is a Lie algebra. When this is the case, the polynomial
algebra part restricts isomorphically to that of $H^*(\Omega_1(G);\mathbb{F}_p)$
and the exterior algebra part is induced isomorphically from that of
$H^*(G/\Omega_1(G);\mathbb{F}_p)$ via the projection homomorphism.
\end{thm}

Fix $G \in \Obj(\mathbf{BGrp})$ we now study the
Lyndon-Hochschild-Serre (L.H.S.) spectral sequence (with
$\mathbb{F}_p$-coefficients) associated
to the extension
$$
1 \rightarrow \Omega_1(G)=V \overset{i}{\rightarrow} G 
\overset{\pi}{\rightarrow} W \rightarrow 1.
$$
Recall that this is a spectral sequence with 
$$
E_2^{p,q}=H^p(W,H^q(V,\mathbb{F}_p))
$$ 
and it abuts to $H^*(G;\mathbb{F}_p)$.
Since this is a central extension, if we use $\mathbb{F}_p$ coefficients 
throughout then 
$$
E_2^{p,q} = H^p(W,H^q(V,\mathbb{F}_p)) 
	  = H^q(V,\mathbb{F}_p) \otimes H^p(W,\mathbb{F}_p) .
$$
So using that $p$ is odd one gets explicitly:
\begin{align}
E_2^{*,*} \cong \wedge^*(V^*) \otimes \mathbb{F}_p[\beta(V^*)]
\otimes \wedge^*(W^*) \otimes \mathbb{F}_p[\beta(W^*)].
\end{align}
As before, we are using basis free notation, for
example $V^*=H^1(V;\mathbb{F}_p)$
and $\beta$ is the Bockstein.

The dual of the $p$-power map $\phi$ gives us an isomorphism
$\phi^* : V^* \rightarrow W^*$ so
given a basis $E=\{ e_1, \dots e_n \}$ for 
$V^*$, one can use $\phi^*(E)=\{x_1,\dots,x_n \}$ as a basis for $W$. 

By standard comparisons it is easy to see
that $d_2|_{W^*}=d_2|_{\beta(W^*)}=d_2|_{\beta(V^*)}=0 $ while
$$
d_2(e_i)=(\text{the ith component of the extension element})=\beta(x_i) 
+ br^*(x_i) 
$$
for all $i$. Here we are using the form of the extension element as presented
in the end of section~\ref{sec: prelim}. 
As $\phi^*(e_i) = x_i$, one has
$$
d_2|_{V^*} = \beta \circ \phi^* + br^* \circ \phi^* .
$$

Let us calculate $E_3^{*,*}$. 

Let 
$$
A = \wedge^*(X_1,\dots,X_n) \otimes \mathbb{F}_p[Y_1,\dots,Y_n] \otimes
\wedge^*(T_1,\dots,T_n) \otimes \mathbb{F}_p[S_1,\dots,S_n]
$$
be an abstract free graded-commutative algebra where the polynomial
generators are degree 2 and the exterior generators are degree 1.
Since this is a free object, the assignment $X_i \rightarrow x_i,
T_i \rightarrow e_i, S_i \rightarrow \beta(e_i), Y_i \rightarrow
\beta(x_i) + br^*(x_i)$ defines a map $\Psi$ of graded-algebras from
$A$ to $E_2^{*,*}$. As $\beta(x_i) + br^*(x_i)$ differs form $\beta(x_i)$
by a nilpotent element it is easy to see by induction on the grading
that $\Psi$ is an isomorphism. The induced differential $D_2$ on $A$ of $d_2$
under this isomorphism has $D_2(X_i)=D_2(Y_i)=D_2(S_i)=0$ and
$D_2(T_i)=Y_i$. Applying K\"unneth's Theorem, we then see that the cohomology
of $(A,D_2)$ is isomorphic to 
$\wedge^*(X_1,\dots,X_n) \otimes \mathbb{F}_p[S_1,\dots,S_n]$.
Using this, we see that $E_3^{*,*}$ is given by
\begin{align}
E_3^{*,*} \cong \wedge^*(W^*) \otimes
\mathbb{F}_p[\beta(V^*)]
\end{align}
We have obviously that $d_3|_{W^*}=0$ and we also have
(see for example page 155 of \cite{Benson}) :
\begin{align}
\begin{split}
\label{eq: one}
d_3 \circ \beta|_{V^*} &= \{ \beta \circ d_2|_{V^*} \} \in 
E_3^{3,0}=\wedge^3(W^*) \\
&= \{ \beta \circ (\beta \circ \phi^* + br^* \circ \phi^*) |_{V^*} \} \\
&= \{ \beta \circ br^* \circ \phi^* |_{V^*} \} \text{ since } \beta \circ \beta = 0.
\end{split}
\end{align}
Now 
note that $\beta$ is a derivation and that $\beta(x_i)$ is identified
with $-br^*(x_i)$ in $E_3^{*,0}=\wedge^*(x_1,\dots,x_n)$ hence $\beta$
induces the same map as $-br^*$ on $E_3^{*,0}$. Thus we see
that the final expression in equation~\ref{eq: one} becomes
$ -br^* \circ br^* \circ \phi^*|_{V^*}$. As $\phi^*$ is an isomorphism,
we can conclude that $d_3|_{\beta(V^*)} = 0$ and hence $d_3 = 0$
if and only if $br^* \circ br^* =0$
on $W^*=\wedge^1(W^*)$ or in other words if and only if $(W^*, br^*)$ is
a co-Lie algebra. Hence we conclude: 

\begin{lem}
$d_3=0$ if and only if
$Log(G) = (W,br)$ is a Lie algebra.
\end{lem}
When this happens $E_3 = E_4$ but then for
grading reasons, $d_r|_{\beta(V^*)}=d_r|_{W^*}=0$ 
for all $r \geq 4$. Hence
as the differentials are derivations which vanish on the generators, all 
further differentials in the spectral sequence must be 0. We have proven
\begin{pro}
$E_3 =E_{\infty}$ if and only if $Log(G)$ is a
Lie algebra. In this case,
$$
E_3^{*,*} =E_{\infty}^{*,*}=\wedge^*(W^*) \otimes
\mathbb{F}_p[\beta(V^
*)].
$$
\end{pro}
Due to the free nature of this graded ring, one can then show in
a routine manner that 
$$
H^*(G;\mathbb{F}_p) = \wedge^*(W^*) \otimes \mathbb{F}_p[S^*]
$$
where we have abused notation a bit and identified $W^*$ with its image
under the map $\pi^*: H^*(W;\mathbb{F}_p) \rightarrow H^*(G;\mathbb{F}_p)$.
$S^*$ is a subspace which 
restricts isomorphically to the subspace $\beta(V^*)$ of 
$H^*(V;\mathbb{F}_p)$. Hence we have proved theorem~\ref{thm: Fp cohomology}.

Let us define
$\mathbf{LGrp}$ to be the full subcategory of $\mathbf{BGrp}$ whose objects are
the bracket groups $G$ where the associated $Log(G)$ is a Lie algebra.
Let $\mathbf{Lie}$ be the full subcategory of $\mathbf{Brak}$ whose objects
are the Lie algebras. Then restricting the natural functors we
had before we see $\mathbf{LGrp}$ and $\mathbf{Lie}$ are naturally equivalent categories.

\subsection{A formula for the Bockstein on $H^*(G ; \mathbb{F}_p)$}
From now on we consider groups $G \in \Obj(\mathbf{LGrp})$. By the previous theorem
$$
H^*(G;\mathbb{F}_p)=\wedge^*(W^*) \otimes \mathbb{F}_p[S^*]
$$
where $W=G/\Omega_1(G).$

We wish to study the Bockstein $\beta : H^*(G;\mathbb{F}_p) \rightarrow H^{*+1}(G,\mathbb{F}_p)$.
The reason for this is that one can determine the groups
$H^*(G;\mathbb{Z})$
from knowledge of $\beta$ and the ``higher Bocksteins'' on 
$H^*(G;\mathbb{F}_p)$.

Recall that $\beta$ is a derivation, i.e.\ , for homogeneous elements
$u,v \in H^*(G;\mathbb{F}_p)$ we have
$$
\beta(uv) = \beta(u)v + (-1)^{\degr(u)}u\beta(v).
$$
Also recall that $\beta \circ \beta = 0$. Since $\beta$ is a derivation
we need only describe it on the generating subspaces $W^*$ and $S^*$
of $H^*(G;\mathbb{F}_p)$. 

Now recall the central short exact sequence
$$
1 \rightarrow \Omega_1(G)=V \overset{i}{\rightarrow} G 
\overset{\pi}{\rightarrow} W \rightarrow 1
$$
Then we have seen that the subalgebra $\wedge^*(W^*)$ of $H^*(G;\mathbb{F}_p)$
is the image of the map 
$\pi^*: H^*(W;\mathbb{F}_p) \rightarrow H^*(G;\mathbb{F}_p)$. 
In the L.H.S.-spectral sequence of the last section, we can identify
this image with $E_{\infty}^{*,0}=E_3^{*,0}$. We saw there that
$\beta$ agrees with $-br^*$ on this subalgebra.

Thus $\beta |_{W^*}= -br^*$. 
It is easy to check that
the differential complex $(\wedge^*(W^*), \beta)$ is the standard 
Koszul resolution used to calculate
the Lie algebra cohomology $H^*(\mathfrak{L};\mathbb{F}_p)$, where
$\mathfrak{L}=(W,br)$. Also, we see that one 
can obtain the Lie algebra structure of $\mathfrak{L}$ from
knowledge of $\beta$ on the exterior algebra part of $H^*(G;\mathbb{F}_p)$. 

Now we
are left with finding $\beta |_{S^*}$ which is harder. Let us develop
some more notation. Purely from considering degree restrictions one has:
$$
\beta |_{S^*} : S^* \rightarrow (W^* \otimes S^*) \oplus \wedge^3(W^*).
$$
Thus we can write 
\begin{equation}
\label{eq: component}
\beta |_{S^*} = \beta_1 + \beta_2 
\end{equation}
where
$\beta_1 : S^* \rightarrow W^* \otimes S^*$ and $\beta_2 : S^* \rightarrow
\wedge^3(W^*)$. If we let $j = \beta^{-1} \circ i^* :S^* \rightarrow
V^* $ be the natural isomorphism, we can define $\bone = (1 \otimes j) \circ
\beta_1 \circ j^{-1} : V^* \rightarrow (W^* \otimes V^*)$ and 
$\btwo = \beta_2 \circ j^{-1} : V^* \rightarrow
\wedge^3(W^*)$. 

	Before we go on, we should note an inherent ambiguity.
Note $\bone$ and $\btwo$ determine $\beta$, however their
definition requires a choice for the subspace $S^*$.  
We would like that all our expressions involving $\beta$
be determined from knowledge of the underlying Lie algebra 
$\mathfrak{L}=Log(G)$.
Note $W=G/\Omega_1(G)$ can be identified with the underlying vector
space of $\mathfrak{L}$ and hence so can $V=\Omega_1(G)$ via the
$p$-power isomorphism $\phi$. Thus the subalgebra $\wedge^*(W^*)$ of 
$H^*(G;\mathbb{F}_p)$ is unambiguously determined by $\mathfrak{L}$.
However the subspace $S^*$ which generates the polynomial part
of $H^*(G;\mathbb{F}_p)$ is only determined
by the fact it restricts isomorphically to the natural subspace
$\beta(V^*)$ of $H^*(V;\mathbb{F}_p)$. We can see that given a
homomorphism $\mu: S^* \rightarrow \wedge^2(W^*)$, then the image
of $S^*$ under $\eta = Id + \mu$ works just as well, and that all possible
other choices for this generating subspace arise in this way. Let
us see how changing ones choices effects the decomposition of $\beta$
in equation~\ref{eq: component}.

It is easiest to see this by using a basis $\{ \bar{s}_1, \dots \bar{s}_n \}$
for the space $S^*$. If we put this basis in a column vector:
$$
\mathbf{s} = \begin{pmatrix} \bar{s}_1 \\ \vdots \\ \bar{s}_n \end{pmatrix}
$$
then we have that
$$
\beta(\mathbf{s}) = \xi \mathbf{s} + \beta_2(\mathbf{s})
$$
where $\xi$ is an $n \times n$ matrix with entries in $\wedge^1(W^*)$.

If we change our choice of subspace $S^*$ using the homomorphism $\mu$
mentioned above, we get a new subspace $\bar{S}^* = \eta(S^*)$ with basis 
$\mathbf{s}'=\mathbf{s} + \mu(\mathbf{s})$.
 We then calculate:
\begin{align*}
\label{eq: alter}
\beta(\mathbf{s}') = \beta(\mathbf{s} + \mu(\mathbf{s})) = 
\xi\mathbf{s} + \beta_2(\mathbf{s}) + \beta \circ \mu(\mathbf{s})
= \xi\mathbf{s}' + (\beta_2(\mathbf{s}) + (-br^* - \xi)(\mu(\mathbf{s}))).
\end{align*}
(Here we have used again that $\beta = -br^*$ on $\wedge^*(W^*)$.)
So if we decompose $\beta = \beta_1' + \beta_2'$ using this new subspace
$\bar{S}^*$ and note that $(1 \wedge \mu) \circ \beta_1(\mathbf{s}) =
\xi \mu(\mathbf{s})$, we see that:
\begin{equation}
\label{eq: etaboundary}
(1 \otimes \eta^{-1}) \circ \beta_1' \circ \eta =\beta_1 \text{ and }
\beta_2' \circ \eta =\beta_2 + (-br^* \circ \mu - (1 \wedge \mu) \circ \beta_1 )
\end{equation}
or
\begin{equation}
\label{eq: etaboundaryII}
\bone' = \bone \text{ and } \btwo' = \btwo + 
(-br^* \circ \mu - (1 \wedge \mu) \circ \bone)
\end{equation}
(Here we have made the natural identification of $\mu$ as a map
$\mu : V^* \rightarrow \wedge^2(W^*)$.)
Thus $\bone$ is not ambiguous 
while $\btwo$ 
depends on the particular subspace $S^*$ chosen, as indicated in 
equation~\ref{eq: etaboundaryII}. We will return to this later. 

Let us now concentrate on $\bone$ which is well defined given
the Lie algebra $\mathfrak{L}$. Since
$\beta \circ \beta = 0$ we have:
\begin{align}
\label{eq: eight}
0 = \beta \circ \beta |_{S^*} = \beta \circ \beta_1 + \beta \circ \beta_2. 
\end{align}
However since $\beta$ is a derivation and $\beta = -br^*$ on $\wedge^*(W^*)$
, it follows easily from equation~\ref{eq: eight} that
\begin{align}
\label{eq: nine}
0 = (-br^* \otimes 1) \circ \beta_1 - (1 \otimes \beta_1) \circ \beta_1
- (1 \otimes \beta_2) \circ \beta_1
-br^* \circ \beta_2.
\end{align}
Taking components in the usual way, we get the equations:
\begin{equation}
\label{eq: coliecomod}
(-br^* \otimes 1) \circ \beta_1 = (1 \otimes \beta_1) \circ \beta_1
\end{equation}
and
\begin{equation}
\label{eq: cohclass}
(1 \otimes \beta_2) \circ \beta_1 = -br^* \circ \beta_2.
\end{equation}
Equation~\ref{eq: coliecomod} gives us the following commutative diagram:
$$
\begin{CD}
V^* @>{-\bone}>> W^* \otimes V^* \\
@V{-\bone}VV @VV{1 \otimes -\bone}V \\
W^* \otimes V^* @>{br^* \otimes 1}>> \wedge^2(W^*) \otimes V^* 
\end{CD}
$$
Thus, $-\bone$ equips $V^*$ with the structure of a comodule 
over the co-Lie algebra $(W^*, br^*)$. Thus, taking duals, we see that
$-\bone^*$ equips $V$ with the structure of a module over the Lie algebra
$(W,br)$. However, we can naturally identify $V$ with $W$ using the
$p$-power map, thus for every Lie algebra $\mathfrak{L} = (W, br)$ we obtain
a map $\lambda: \mathfrak{L} \rightarrow \mathfrak{gl}\hskip.02in
(\mathfrak{L})$ from the module structure above. Furthermore it follows,
that if we
give $\mathfrak{gl}\hskip.02in(\mathfrak{L})$ the canonical Lie
algebra structure, then this map
$\lambda$ will be a map of Lie algebras.
 
These maps $\lambda$ obtained above are natural in the following sense:
\begin{claim}
\label{claim: naturality}
Let $\mathfrak{L},\mathfrak{L}' \in \Obj(\mathbf{Lie})$ and $\psi \in 
Mor(\mathfrak{L},\mathfrak{L}')$ then 
$$
[\lambda'(\psi(x))](\psi(y)) = \psi([\lambda(x)](y)).
$$
(The notation $[\lambda(x)](y)$ means of
course $\lambda(x) \in \mathfrak{gl}\hskip.02in(\mathfrak{L})$ applied
to the element $y$.)
\end{claim}
\begin{proof}
Suppose $\mathfrak{L}=(W,br)$ and $\mathfrak{L}'=(\bar{W},\bar{br})$. Then
$\psi: W \rightarrow \bar{W}$ induces a unique morphism $Exp(\psi) :
Exp(\mathfrak{L}) \rightarrow Exp(\mathfrak{L}')$ in the category
$\mathbf{LGrp}$. This is a class of group
homomorphisms where any two representatives of this class induce the same
mapping on the $\Omega_1$ level as mentioned earlier. It is easy
to check that any representative of this class will induce the 
dual map $\psi^*$ from 
$$
\bar{W}^*=H^1(Exp(\mathfrak{L}');\mathbb{F}_p) \rightarrow
W^*=H^1(Exp(\mathfrak{L}); \mathbb{F}_p).
$$
Since the Bockstein is a natural operation, it will
commute with any such map on cohomology and from this, one obtains the
following commutative diagram:
$$
\begin{CD}
\bar{V}^* @>{-\bar{\bone}}>> \bar{W}^* \otimes \bar{V}^* \\
@V{\psi^*}VV @VV{\psi^* \otimes \psi^*}V \\
V^* @>{-\bone}>> W^* \otimes V^*
\end{CD}
$$
The claim follows easily now by taking the dual of the diagram above
and recalling that the map $\lambda$ was defined via $-\bone^*$. 
\end{proof}

Notice that knowledge of the map $\lambda$ is equivalent to knowledge
of the component $\bone$ of the Bockstein.
We have isolated enough properties of the maps $\lambda$ and will
now set out to determine them explicitly.

\subsection{Self-representations}
In this section, let $\mathbf{k}$ be a finite field with 
$\chan(\mathbf{k}) \neq 2$. Let $\mathbf{Lie}$ be the
category of finite dimensional Lie algebras over $\mathbf{k}$ with morphisms
the lie algebra maps.
(We will of course be interested in the case $k=\mathbb{F}_p$).
\smallskip
\begin{defn}
A natural self-representation on $\mathbf{Lie}$ is a collection of linear maps
$\kappa_{\mathfrak{L}} : \mathfrak{L} \rightarrow \mathfrak{gl}\hskip.02in(\mathfrak{L})$, which
satisfy the naturality condition
$$
[\kappa_{\mathfrak{L}'}(\psi(x))](\psi(y)) = \psi([\kappa_{\mathfrak{L}}(x)](y))
$$
for all $\mathfrak{L},\mathfrak{L}' \in \Obj(\mathbf{Lie})$,  
$\psi \in \Mor(\mathfrak{L},\mathfrak{L}')$ and $x,y \in \mathfrak{L}$.
One says that a self-representation is strong if the maps
$\kappa_{\mathfrak{L}} : \mathfrak{L} \rightarrow
\mathfrak{gl}\hskip.02in(\mathfrak{L})$
are maps of Lie algebras.
\end{defn}

\noindent
\textbf{Examples:} 

\noindent
(a) If we set $\kappa = 0$ we see easily that we get a strong self-representation
which we call the zero representation. 
\smallbreak
\noindent
(b) By claim~\ref{claim: naturality},
we see that the maps $\lambda$ fit together to give a
strong natural self-representation which we will call $\lambda$.
\smallbreak
\noindent 
(c) If we set $\kappa = ad$ where 
$ad : \mathfrak{L} \rightarrow \mathfrak{gl}\hskip.02in(\mathfrak{L})$ is given by
$ad(x)(y) = [x,y]$ for all $x,y \in \mathfrak{L}$, then it is easy to see
that each map ad is a map of Lie algebras. Furthermore if $\mathfrak{L}'$
is another Lie algebra and $\psi : \mathfrak{L} \rightarrow \mathfrak{L'}$
is a morphism of Lie algebras then we have for $x,y \in \mathfrak{L}$ : 
$$
[ad'(\psi(x))](\psi(y)) = [\psi(x),\psi(y)] = \psi([x,y]) = \psi([ad(x)](y)).
$$
So we see this assignment defines a strong natural self-representation
which we call the adjoint representation.

Now we will study self-representations so that we can show $\lambda$
is the adjoint representation. Note that if $\kappa_0,\kappa_1$ are two
self-representations, then if they assign the same map
$\mathfrak{L} \rightarrow \mathfrak{gl}\hskip.02in(\mathfrak{L})$ for some lie algebra
$\mathfrak{L}$, (we will say they agree on $\mathfrak{L}$) then they assign the 
same map for all lie algebras isomorphic to $\mathfrak{L}$. So we will implicitly 
use this from now on without mention.

\begin{lem}
\label{lem: subalgebra}
If $\kappa_0,\kappa_1$ are two self-representations, then if
they agree on a lie algebra $\mathfrak{L}$, they agree on all 
Lie subalgebras of $\mathfrak{L}$.
\end{lem}
\begin{proof} 
Follows from the naturality of self-representations.
\end{proof}

Recall $\mathfrak{gl}_n$ is the Lie algebra of $n \times n$ matrices
and $\mathfrak{sl}_n$ is the lie algebra of $n \times n$ matrices with
trace zero. Define $\mathfrak{N}$ to be the 3-dimensional Lie algebra
with basis $\{x,y,z\}$ and bracket given uniquely by $[x,y]=z,[x,z]=[y,z]=0$. 
Note that $x,z$ generate a 2-dimensional abelian subalgebra of $\mathfrak{N}$.
Define $\mathfrak{S}$ to be the 2-dimensional Lie algebra with basis
$\{x,y\}$ and bracket given uniquely by $[x,y]=x$. Note as 
$\chan(\mathbf{k}) \neq 2$
this Lie algebra embeds into $\mathfrak{sl}_2$ by the Lie algebra map $\Psi$ where
$$
\Psi(x) = \begin{pmatrix}  0 & 0 \\ 1 & 0 \end{pmatrix}  
$$
$$
\Psi(y) = \begin{pmatrix}  1/2 & 0 \\ 0 & -1/2 \end{pmatrix}.
$$

\begin{lem}
\label{lem: reduction}
Let $\kappa_0,\kappa_1$ be two self-representations, then if
$\kappa_0$ and $\kappa_1$ agree on $\mathfrak{sl}_2$ and $\mathfrak{N}$
then $\kappa_0=\kappa_1$ as self-representations.
\end{lem}
\begin{proof}
Assume $\kappa_1,\kappa_0$ are two self-representations satisfying
the assumptions above.
By the Ado-Iwasawa theorem, all Lie algebras embed as subalgebras of
$\mathfrak{gl}_n$ for some $n$. So by lemma~\ref{lem: subalgebra}, 
to show $\kappa_1=\kappa_0$
one need only show they agree on $\mathfrak{gl}_n$ for all $n$.
A basis for $\mathfrak{gl}_n$ is $\{ \delta_{i,j} : i,j =1,\dots n \}$ where
$\delta_{i,j}$ is the matrix with entry one at the (i,j)-position and
zero elsewhere. One calculates easily that:
$$
[\delta_{i,j},\delta_{l,m}] =  
\begin{cases}
0 &\text{ if } j \neq l, m \neq i \\
\delta_{i,m} &\text{ if } j=l, m \neq i \\
-\delta_{l,j} &\text{ if } j \neq l, m=i \\
\delta_{i,i} - \delta_{j,j} &\text{ if } j=l,m=i.
\end{cases}
$$
Let $\mathfrak{L}_{ij,lm}^n$ be the Lie subalgebra of $\mathfrak{gl}_n$ generated 
by $\delta_{i,j}$ and $\delta_{l,m}$ for all $i,j,l,m \in \{1,\dots, n \}$.
Suppose we have shown $\kappa_1$ and $\kappa_0$ agree on $\mathfrak{L}_{ij,lm}^n$
for all $i,j,l,m \in \{1,\dots,n \}$. Then by an easy linearity argument
it follows that $\kappa_1$ and $\kappa_0$ agree on $\mathfrak{gl}_n$ for all $n$ and hence are equal. So we see that if we can show $\kappa_1$ and $\kappa_0$ 
agree on $\mathfrak{L}_{ij,lm}^n$ for all 
$n \in \mathbb{N}, i,j,l,m \in \{1,\dots,n\}$ we are done.
(Since $\mathfrak{gl}_n \subset \mathfrak{gl}_{n+1}$ we see we can assume $n \geq 5$
say.) Now let us identify the Lie algebras $\mathfrak{L}_{ij,lm}^n$. We have
the following cases:
\begin{enumerate}
\item $\mathfrak{L}_{ij,lm}$ where $j \neq l, m \neq i$.
Here we have $[\delta_{i,j},\delta_{l,m}]=0$ so this Lie algebra is
an abelian Lie algebra of dimension $\leq 2$. Hence it is contained
inside $\mathfrak{N}$.
\item $\mathfrak{L}_{ij,jm}$ where $m \neq i, m \neq j, j \neq i$.
This Lie algebra has basis $\{ \delta_{i,j},
\delta_{j,m}, \delta_{i,m} \}$  and has bracket relations,
$[\delta_{i,j},\delta_{j,m}]=\delta_{i,m}$, $[\delta_{i,j},\delta_{i,m}]=
[\delta_{j,m},\delta_{i,m}]=0$. Thus it is easy to see this Lie algebra is
isomorphic to $\mathfrak{N}$.
\item $\mathfrak{L}_{ij,jj}$ where $j \neq i$.
Here we have $[\delta_{i,j},\delta_{j,j}]=\delta_{i,j}$. So we see
this Lie algebra is isomorphic to $\mathfrak{S}$ and so is inside $\mathfrak{sl}_2$.
\item $\mathfrak{L}_{jj,jm}$ where $m \neq j$.
We see easily we get the same case as in case 3.
\item $\mathfrak{L}_{ij,li}$ where $j \neq l$.
Here we note that $\mathfrak{L}_{ij,li} = \mathfrak{L}_{li,ij}$ so it occurs
as one of the cases 2 -- 4.
\item $\mathfrak{L}_{ij,ji}$ where $i \neq j$.
Here we have that $\mathfrak{L}_{ij,ji}$ has basis $\{ \delta_{i,j}, \delta_{j,i},
\delta_{i,i}-\delta_{j,j}=\Delta \}$. The bracket is given by
$[\delta_{i,j},\delta_{j,i}]=\Delta$, $[\Delta,\delta_{i,j}]
=\delta_{i,j}+\delta_{i,j}=2\delta_{i,j}$ and $[\Delta,\delta_{j,i}]
=-\delta_{j,i}-\delta_{j,i}=-2\delta_{j,i}$ and so we see easily
that this Lie algebra is isomorphic to $\mathfrak{sl}_2$.
\item $\mathfrak{L}_{ii,ii}$.
Obviously this Lie algebra is 1-dimensional and hence lies in either
$\mathfrak{sl}_2$ or $\mathfrak{N}$.
\end{enumerate}

\smallbreak
\noindent
All the Lie algebras $\mathfrak{L}_{ij,lm}^n$ fit into one of these cases
and hence embed in either $\mathfrak{sl}_2$
or $\mathfrak{N}$. As $\kappa_1$ and $\kappa_0$ agree on these two
Lie algebras by assumption, they must agree on all the 
$\mathfrak{L}_{ij,lm}^n$ and hence
must be equal as argued before.
\end{proof}

Now we will study how a self-representation must look in certain particular
cases.
\begin{lem}
\label{lem: abelian}
Let $\kappa$ be a self-representation. If $x,y \in \mathfrak{L}$ have $[x,y]=0$ then $[\kappa(x)](y)=0$. So in particular $\kappa=0$ on abelian Lie algebras.
\end{lem}
\begin{proof}
By naturality, it is enough to show $\kappa=0$ on abelian Lie algebras
of dimension $\geq 2$. Fix $\mathfrak{L}$ an abelian Lie algebra of
dimension $n \geq 2$.
Let $\gamma \in \mathbf{k}-\{0,1\}$.
Then multiplication by $\gamma$ induces an automorphism $\Gamma$ 
of the Lie algebra
$\mathfrak{L}$. Given $x, y \in \mathfrak{L}$ then if
$[\kappa(x)](y)=z$ it follows from naturality that 
\begin{align*}
\begin{split}
[\kappa(\Gamma(x))](\Gamma(y)) &= \Gamma(z) \\
\gamma^2 z = \gamma z 
\end{split}
\end{align*}
However, $\gamma^2 \neq \gamma$ so $z=[\kappa(x)](y)=0$. Since $x,y$ were
arbitrary, the lemma follows.
\end{proof}
\begin{cor}
\label{cor: basic}
Let $\kappa$ be a self-representation.
Then for a Lie algebra $\mathfrak{L}$ and $x,y \in \mathfrak{L}$ one has
$[\kappa(x)](x) =0$ and $[\kappa(x)](y)= -[\kappa(y)](x)$.
\end{cor}
\begin{proof}
As $[x,x]=0$, $[\kappa(x)](x)=0$ follows from lemma~\ref{lem: abelian}.
The second part follows easily from $[\kappa(x+y)](x+y) = 0$. 
\end{proof}  

Recall
$\mathfrak{S}$ is the Lie algebra with basis $\{x,y\}$ and $[x,y]=x$.
Let $\pi : \mathfrak{S} \rightarrow \mathfrak{S}$ be a linear map given by
$\pi(x)= 0$ and $\pi(y)= y$. It is easy to check that $\pi$ is a map of Lie 
algebras.

\begin{lem}
\label{lem: S}
Let $\kappa$ be a self-representation.
Let $\mathfrak{S}$ be the Lie algebra with basis $\{x,y\}$ and bracket
given by $[x,y]=x$. Then there is a number $a(\kappa) \in \mathbf{k}$ such
that one has 
$[\kappa(u)](v) = a(\kappa)[u,v]$ for all $u,v \in \mathfrak{S}$. 
\end{lem}
\begin{proof}
By corollary~\ref{cor: basic}, one sees it is enough to show
$[\kappa(x)](y)=a(\kappa)x$ for some $a(\kappa) \in \mathbf{k}$.
One knows a~priori that
\begin{align}
\label{eq: kappa1}
[\kappa(x)](y) = a(\kappa)x + b(\kappa)y,
\end{align}
so we want to show $b(\kappa)=0$. Applying the map of Lie algebras
$\pi$, mentioned in the paragraph preceding the lemma,
to equation~\ref{eq: kappa1} and using naturality, we get
$0 = b(\kappa)y$ and so it follows that $b(\kappa)=0$.
\end{proof}
Recall $\mathfrak{sl}_2$ has basis $\{ \delta_{1,2},\delta_{2,1},
\delta_{1,1}-\delta_{2,2}\}$. Let us relabel as follows:
$H=\delta_{1,1}-\delta_{2,2}$, $X_{+} = \delta_{1,2}$ and
$X_{-} = \delta_{2,1}$ then it is easy to verify that the bracket on
$\mathfrak{sl}_2$ is given by 
\begin{align}
[X_{+},X_{-}] = H ,
[H,X_{+}] = 2X_{+} \text{ and } 
[H,X_{-}] = -2X_{-}.
\end{align}
Define linear maps $\Psi_{mix} : \mathfrak{sl}_2 \rightarrow
\mathfrak{sl}_2$, $\Psi_{neg} : \mathfrak{sl}_2 \rightarrow \mathfrak{sl}_2$ as
follows: 
\begin{align}
\Psi_{neg}(H) = H , 
\Psi_{neg}(X_{+}) = -X_{+} ,
\Psi_{neg}(X_{-}) = -X_{-} , 
\end{align}
and
\begin{align}
\Psi_{mix}(H) = H + 2X_{+} , 
\Psi_{mix}(X_{+}) = -X_{+} ,
\Psi_{mix}(X_{-}) = H + X_{+} -X_{-}.
\end{align}
It is easy to see that $\Psi_{neg}$ is a map of Lie algebras, and a routine
calculation shows $\Psi_{mix}$ is one as well.
Now we are ready to prove:
\begin{lem}
\label{lem: sl2}
Let $\kappa$ be a self-representation. Then there is $a(\kappa) \in \mathbf{k}$ which is the same as in
lemma~\ref{lem: S}, such that
$\kappa=a(\kappa) ad$ on $\mathfrak{sl}_2$, i.e., 
$[\kappa(x)](y)=a(\kappa)[x,y]$
for all $x,y \in \mathfrak{sl}_2$. Furthermore, if $\kappa$ is strong, then 
$a(\kappa)=0$ or $1$.
\end{lem}
\begin{proof}
One has $[\kappa(H)](H)=[\kappa(X_{+})](X_{+})=[\kappa(X_{-})](X_{-})=0$
as usual. Now $H, X_{+}$ span a Lie subalgebra of $\mathfrak{sl}_2$ which
is isomorphic to $\mathfrak{S}$ via $H \rightarrow -2y$ and $X_{+} \rightarrow
x$. Thus by naturality and lemma~\ref{lem: S} we have 
$[\kappa(H)](X_{+})=a(\kappa)[H,X_{+}]$. A similar argument works for $X_{-}$
in place for $X_{+}$.
So we conclude
easily that
\begin{align}
\label{eq: kappaeq}
[\kappa(H)](\cdot) = a(\kappa)[H,\cdot]. 
\end{align}
Thus it follows that
$[\kappa(X_{\pm})](H)=a(\kappa)[X_{\pm},H]$.
Now $[\kappa(X_{+})](X_{-}) = 
\alpha H + \beta X_{+} + \gamma X_{-}$. Applying the map of Lie
algebras $\Psi_{neg}$ to this equation, and using naturality one gets:
\begin{align*}
\begin{split}
[\kappa(\Psi_{neg}(X_{+}))](\Psi_{neg}(X_{-})) =&
\alpha H - \beta X_{+} - \gamma X_{-} \\
\alpha H + \beta X_{+} + \gamma X_{-} =& \alpha H - \beta X_{+}
- \gamma X_{-}.
\end{split}
\end{align*}
So $\beta = \gamma=0$. Thus $[\kappa(X_{+})](X_{-})=\alpha H =\alpha 
[X_{+},X_{-}]$. So if we can show $\alpha = a(\kappa)$ we can conclude
that
$$
[\kappa(X_{\pm})](\cdot) = a(\kappa)[X_{\pm},\cdot]
$$
and together with equation~\ref{eq: kappaeq} this implies
 $\kappa=a(\kappa) ad$ on 
$\mathfrak{sl}_2$. So it remains only to show $\alpha = a(\kappa)$.
Applying the Lie algebra map $\Psi_{mix}$ to the equation
$[\kappa(X_{+})](X_{-})=\alpha H$ one gets:
\begin{align*}
\begin{split}
[\kappa(\Psi_{mix}(X_{+}))](\Psi_{mix}(X_{-})) =&
\alpha \Psi_{mix}(H) \\
-[\kappa(X_{+})](H) + [\kappa(X_{+})](X_{-}) =&
\alpha (H + 2X_{+}) \\
-a(\kappa)[X_{+},H] + \alpha H =& \alpha H + 2\alpha X_{+} \\
2a(\kappa)X_{+} =& 2\alpha X_{+}.
\end{split}
\end{align*}
From which one concludes that $a(\kappa)=\alpha$ as desired. 

Now if
$\kappa=a(\kappa)ad$ is a map of Lie algebras, then we have
\begin{align*}
\begin{split}
\kappa([X_+,X_-]) &= \kappa(X_+) \circ \kappa(X_-) - \kappa(X_-) \circ
\kappa(X_+) \\
a(\kappa)ad([X_+,X_-]) &= a(\kappa)^2 ad(X_+) \circ ad(X_-) - a(\kappa)^2 
ad(X_-)
\circ ad(X_+) \\
a(\kappa)ad([X_+,X_-]) &= a(\kappa)^2 ad([X_+,X_-]).
\end{split}
\end{align*}
Since $ad([X_+,X_-])$ is nonzero, it follows that $a(\kappa) = a(\kappa)^2$
or in other words, that $a(\kappa)=0$ or $1$.
\end{proof}

Recall $\mathfrak{N}$ has basis $\{x,y,z\}$ and bracket given by
$[x,z]=[y,z]=0, [x,y]=z$. Define linear maps 
$\Psi_x,\Psi_y : \mathfrak{N} \rightarrow \mathfrak{N}$ by 
$$
\Psi_x(y) = \Psi_x(z) = 0, \Psi_x(x) = x
$$ 
and 
$$
\Psi_y(x) = \Psi_y(z) = 0, \Psi_y(y) = y.
$$
It is easy to verify that these maps are maps of Lie algebras.
\begin{lem}
\label{lem: N}
Let $\kappa$ be a self-representation.
Then there is $\mu(\kappa) \in \mathbf{k}$ such that
$\kappa = \mu(\kappa) ad$ on $\mathfrak{N}$.
\end{lem}
\begin{proof}
As $z$ is central one sees by lemma~\ref{lem: abelian}, that
$\kappa(z)=0=ad(z)$. Now $[\kappa(x)](y)= \alpha x + \beta y + \mu z$.
Applying the Lie algebra map $\Psi_x$ to this equation and using
naturality one gets $0 = \alpha x$.
So one concludes $\alpha=0$. Similarly using $\Psi_y$ instead one concludes
$\beta = 0$. So $[\kappa(x)](y)=\mu z=\mu [x,y]$. So now it is easy
to see $[\kappa(x)](\cdot)=\mu[x,\cdot]$ and $[\kappa(y)](\cdot)=\mu[y,\cdot]$
via lemma~\ref{lem: abelian} and corollary~\ref{cor: basic}. Thus
setting $\mu(\kappa)=\mu$ we are done.
\end{proof}
Recall $\mathfrak{so}_n \subset \mathfrak{gl}_n$, the subset of skew-symmetric
$n \times n$ matrices, is a Lie subalgebra. For $\chan(k) \neq 2$,
$\mathfrak{so}_3$ is 3-dimensional with basis $\{X=\delta_{1,2}-\delta_{2,1},
Y=\delta_{1,3}-\delta_{3,1},Z=\delta_{2,3}-\delta_{3,2}\}$ and bracket
given by the relations $[X,Y]=-Z$,$[Y,Z]=-X$,$[Z,X]=-Y$. As we are over
a finite field of odd characteristic, $\mathfrak{so}_3$ is isomorphic to
$\mathfrak{sl}_2$ via the map which sends $X \rightarrow \frac{X_+ - X_-}{2},
Y \rightarrow a \frac{X_+ + X_-}{2} + b \frac{H}{2}, Z \rightarrow
b \frac{X_+ + X_-}{2} - a \frac{H}{2}$. Here $a, b \in \mathbf{k}$ are
such that $a^2 + b^2 = -1$. 
\begin{lem}
\label{lem: so3}
Let $\kappa$ be a self-representation.
Then $\kappa=\mu(\kappa) ad=a(\kappa) ad$ on $\mathfrak{so}_3$ where 
$\mu(\kappa)$ is the
same as that of lemma~\ref{lem: N} and $a(\kappa)$ is the same as
lemma~\ref{lem: sl2}. In particular, $a(\kappa) = \mu(\kappa)$.
\end{lem}
\begin{proof}
The fact that $\kappa=a(\kappa) ad$ follows from lemma~\ref{lem: sl2}
and the fact that $\mathfrak{so}_3$ is isomorphic to $\mathfrak{sl}_2$.
So it remains to show that $\kappa=\mu(\kappa) ad$ also.
Consider the Lie algebra $\mathfrak{gl}_3$. In the proof of 
lemma~\ref{lem: reduction} we saw that 
$$
\{ \delta_{i,j},\delta_{j,k},\delta_{i,k} \}
$$
formed a basis of a Lie subalgebra isomorphic to $\mathfrak{N}$ for $i,j,k$
distinct.
Thus 
$$
[\kappa(\delta_{i,j})](\delta_{j,k})=
\mu(\kappa)[\delta_{i,j},\delta_{j,k}]=\mu(\kappa)\delta_{i,k}
$$ 
for all $i,j,k$ distinct. Also $[\kappa(\delta_{i,j})](\delta_{l,m})=0$
if $j \neq l,i \neq m$ as under these conditions 
$[\delta_{i,j},\delta_{l,m}]=0$. So
\begin{align*}
\begin{split}
[\kappa(\delta_{1,2}-\delta_{2,1})](\delta_{1,3}-\delta_{3,1}) =&
[\kappa(\delta_{1,2})](\delta_{1,3}) - [\kappa(\delta_{1,2})](\delta_{3,1}) \\
&-[\kappa(\delta_{2,1})](\delta_{1,3}) +[\kappa(\delta_{2,1})](\delta_{3,1}) \\
=& 0 - (-\mu(\kappa)\delta_{3,2}) -(\mu(\kappa)\delta_{2,3}) + 0 \\
=&\mu(\kappa)(\delta_{3,2}-\delta_{2,3}) \\
=&\mu(\kappa)[\delta_{1,2}-\delta_{2,1},\delta_{1,3}-\delta_{3,1}],
\end{split}
\end{align*}
and similarly
$$
[\kappa(\delta_{1,2}-\delta_{2,1})](\delta_{2,3}-\delta_{3,2}) =
\mu(\kappa)[\delta_{1,2}-\delta_{2,1},\delta_{2,3}-\delta_{3,2}] 
$$
and
$$
[\kappa(\delta_{1,3}-\delta_{3,1})](\delta_{2,3}-\delta_{3,2}) =
\mu(\kappa)[\delta_{1,3}-\delta_{3,1},\delta_{2,3}-\delta_{3,2}].
$$
So with these 3 equations and corollary~\ref{cor: basic} one concludes
that $\kappa =\mu(\kappa) ad$ on $\mathfrak{so}_3$.
\end{proof}

We now state:
\begin{thm}
\label{thm: ad or 0}
Let $\kappa$ be a strong self-representation. 
If $\kappa \neq 0$ on $\mathfrak{sl}_2$ then  
$\kappa = ad$ as self-representations. If $\kappa = 0$ on 
$\mathfrak{sl}_2$ then $\kappa = 0$ as self-representations.
\end{thm}
\begin{proof}
If $\kappa \neq 0$ on $\mathfrak{sl}_2$ then 
we have $a(\kappa)=\mu(\kappa)=1$ from lemmas~\ref{lem: sl2} and 
\ref{lem: so3}. 
Thus $\kappa = ad$ on $\mathfrak{sl}_2$ and on
$\mathfrak{N}$. Thus by lemma~\ref{lem: reduction} we have that
$\kappa=ad$ as self-representations. the case $\kappa=0$ on $\mathfrak{sl}_2$
proceeds similarly.
\end{proof}

For use later on, we need to state theorem~\ref{thm: ad or
0}
in slightly more generality. To do this we need to introduce the 
categories $\mathbf{Lie}_k(p)$ for $k \geq 1$ and $p$ a prime. 
The objects of this category will
be Lie algebras over $\mathbb{Z}/p^k\mathbb{Z}$, that is
free $\mathbb{Z}/p^k\mathbb{Z}$-modules $R$ of finite rank equipped with
a bilinear, alternating form $[\cdot,\cdot] : R \times R \rightarrow R$
which satisfies the Jacobi identity. The morphisms will be the
$\mathbb{Z}/p^k\mathbb{Z}$-module maps which preserve the brackets.
Thus $\mathbf{Lie}_1(p)$ is just the same category $\mathbf{Lie}$
considered before for $\mathbf{k}=\mathbb{F}_p$. 

For $t \leq k$,
we have reduction functors $\mathbf{Lie}_k(p) \rightarrow
\mathbf{Lie}_t(p)$.
These reduction functors are obtained in the obvious way from the reduction map
of $\mathbb{Z}/p^k\mathbb{Z}$ to $\mathbb{Z}/p^t\mathbb{Z}$.
For $\mathfrak{L} \in \Obj(\mathbf{Lie}_k(p))$ we denote
$\bar{\mathfrak{L}}$ to be its reduction in $\Obj(\mathbf{Lie}_1(p))$.
We will also use the bar notation for the reduction of a morphism.
Note that there is the canonical reduction homomorphism which is a map
of abelian groups which preserves
the bracket from $\mathfrak{L}$ to $\bar{\mathfrak{L}}$, we will refer
to maps of abelian groups which preserves brackets as Lie algebra
maps for the rest of this section.

\begin{defn}
A natural self-representation on $\mathbf{Lie}_k(p)$ is a collection of
maps $\kappa_{\mathfrak{L}} : \mathfrak{L} \rightarrow
\mathfrak{gl}(\bar{\mathfrak{L}})$ of abelian groups, 
which satisfy the naturality condition
$$
[\kappa_{\mathfrak{L}'}(\psi(x))](\bar{\psi}(y))
= \bar{\psi}([\kappa_{\mathfrak{L}}(x)](y))
$$
for all $\mathfrak{L}, \mathfrak{L}' \in \Obj(\mathbf{Lie}_k(p))$,
$\psi \in \Mor(\mathfrak{L},\mathfrak{L}')$, $x \in \mathfrak{L}$
and $y \in \bar{\mathfrak{L}}$. We say the natural self-representation
is strong if the maps $\kappa_{\mathfrak{L}} : \mathfrak{L}
\rightarrow \mathfrak{gl}(\bar{\mathfrak{L}})$ are maps of Lie algebras.
\end{defn}
Note this definition extends the previous one. 

It is easy to check
that if we assign $\kappa_{\mathfrak{L}} : \mathfrak{L} \rightarrow
\mathfrak{gl}(\bar{\mathfrak{L}})$ to be the composition of the
reduction homomorphism with $ad : \bar{\mathfrak{L}} \rightarrow
\mathfrak{gl}(\bar{\mathfrak{L}})$ then we obtain a strong
self-representation which we will call $ad$.

The Lie algebras $\mathfrak{gl}_n$, 
$\mathfrak{sl}_2$,$\mathfrak{N}, \mathfrak{S}$ and $\mathfrak{so}_3$
will denote the obvious Lie algebras in $\Obj(\mathbf{Lie}_k(p))$.
For example, $\mathfrak{gl}_n$ is the Lie algebra of $n \times n$
matrices with $\mathbb{Z}/p^k\mathbb{Z}$ entries. 
Theorem~\ref{thm: ad or 0} is true for these
more general self-representations. We state this in the next theorem.

\begin{thm}
\label{thm: ad or 0 II}
Let $p$ be an odd prime and let $\kappa$ be a strong
self-representation over $\mathbf{Lie}_k(p)$. If $\kappa \neq 0$ on
$\mathfrak{so}_3$ and $\mathfrak{sl}_2$, then $\kappa = ad$ as
self-representations. 
\end{thm}
\begin{proof}
The proof proceeds as before with only some 
modifications so we leave it to the reader (see \cite{pak}). 
The only thing to note is
there is a generalized
Ado-Iwasawa theorem (see \cite{W4}) which says that any $\mathfrak{L} \in 
\Obj(\mathbf{Lie}_k(p))$ embeds as a subalgebra of 
$\mathfrak{gl}_n(\mathbb{Z}/p^k\mathbb{Z})$.

\end{proof}

\subsection{Bockstein spectral sequence implies $\lambda=ad$}

Recall we had $\lambda$ a strong self-representation over
$\mathbb{F}_p$
which was defined using the Bockstein. In this section we show 
$\lambda=ad$ as self-representations. 

By Theorem~\ref{thm: ad or 0},
it is enough to show $\lambda$ is not identically zero on $\mathfrak{sl}_2$.
To do this we will need the Bockstein spectral-sequence.

Let X be a topological space (with finitely generated integral homology).
Then the coefficient sequence
$$
0 \rightarrow \mathbb{Z} \overset{p}{\rightarrow} \mathbb{Z} \rightarrow 
\mathbb{Z}/p\mathbb{Z} \rightarrow 0
$$
gives a long exact sequence 
$$ 
\cdots \longrightarrow H^*(X;\mathbb{Z}) \overset{p}{\longrightarrow} 
H^*(X;\mathbb{Z}) 
\overset{mod}{\longrightarrow} H^*(X;\mathbb{F}_p)
\overset{\delta}{\longrightarrow} H^{*+1}(X,\mathbb{Z}) \longrightarrow \cdots
$$
which gives a spectral sequence in a standard way if we think of it as
an exact couple. This spectral sequence
has $B_1^*=H^*(X;\mathbb{F}_p)$ with differential given by the Bockstein
$\beta$ (as defined previously) and converges to $\mathbb{F}_p \otimes
(H^*(X;\mathbb{Z})/\text{torsion})$. If we apply this to the case X=BG where
G is a finite group, then the Bockstein spectral sequence must converge to
zero in positive gradings, as the integral cohomology of a finite group
is torsion in positive gradings. 

Now we will apply the Bockstein spectral sequence 
to show that $\lambda$ is not zero
on $\mathfrak{sl}_2$.

Let $\{H,X_{+},X_{-}\}$ be the usual basis for $\mathfrak{sl}_2$ and
let $G=Exp(\mathfrak{sl}_2)$. Let $\{h,x_{+},x_{-}\}$ be the obvious
dual basis for $\mathfrak{sl}_2$. Then as we have seen before,
$$
B_1^*=H^*(G;\mathbb{F}_p)=\wedge^*(h,x_{+},x_{-}) \otimes
\mathbb{F}_p[s_h,s_+,s_-]
$$
where we have chosen the $s$-basis to correspond in the obvious way.
Then we see using the formulas found before that on the exterior algebra part 
the Bockstein is given by
\begin{align*}
\begin{split}
\beta(h) &= -[\cdot,\cdot]_H=-x_+x_- \\
\beta(x_+) &= -[\cdot,\cdot]_{X_+}=-2hx_+ \\
\beta(x_-) &= -[\cdot,\cdot]_{X_-}=2hx_-.
\end{split}
\end{align*}
So as $B_1^1=\wedge^1(h,x_+,x_-)$ one sees that $B_2^1=0$.
Now note $B_1^2=\wedge^2(h,x_+,x_-) \oplus \spant(s_h,s_+,s_-)$ and by the above
formulas, the exterior part contributes nothing to $B_2^2$. 
If $\lambda$ is zero
on $\mathfrak{sl}_2$ then this means
$\beta$ maps $\spant(s_h,s_+,s_-)$ into 
 $\wedge^3(h,x_+,x_-)$ which is 1-dimensional generated by
$hx_+x_-$. So we see $\beta$ must 
have at least a 2-dimensional kernel on the span of $\{s_h,s_+,s_-\}$.
Thus $\dime(B_2^2) \geq 2$. Now $B_1^3$ is equal to
$$
\wedge^3(h,x_+,x_-) \oplus \spant(\{s_hh,s_hx_+,s_hx_-,s_+h,s_+x_+,
s_+x_-,s_-h,s_-x_+,s_-x_-\}).
$$ 
We will show $\beta$ is injective on
the second part so that $\dime(B_2^3) \leq 1$. Note that
$\beta(s_*)=a_*hx_+x_-$, for some $a_* \in \mathbb{F}_p$, where * stands for 
an arbitrary subscript. So $\beta(s_*)x_*=0$ and one has:
\begin{align*}
\begin{split}
\beta(s_*x_*)=& \beta(s_*)x_* + s_*\beta(x_*) \\ 
=& s_*\beta(x_*) 
\end{split}
\end{align*}
So because $\{\beta(x_h),\beta(x_+),\beta(x_-)\}$ are linearly independent,
and $\{s_h,s_+,s_-\}$ are algebraically independent, it is easy to
see $\beta$ is injective on the second part of $B_1^3$ as claimed.
Thus we have $B_2^1=0$, $\dime(B_2^2) \geq 2$ and $\dime(B_2^3) \leq 1$.
However the Bockstein spectral sequence converges to zero in positive
gradings implying we need $\dime(B_2^2) \leq \dime(B_2^3)$ which we don't
have. Thus we have a contradiction and our assumption that $\lambda=0$
must be false.

Thus we finally conclude $\lambda = ad$ for
all Lie algebras. Now let us chase through the definitions to get formulas
for the $\beta_1$ part of $\beta$. So fix a Lie algebra $\mathfrak{L}$ over
$\mathbb{F}_p$. Fix a basis $E=\{e_1,\dots,e_n\}$ of $\mathfrak{L}$ 
and let $G=Exp(\mathfrak{L})$.
Then $H^*(G;\mathbb{F}_p)=\wedge^*(x_1,\dots,x_n) \otimes \mathbb{F}_p[s_1,
\dots,s_n]$ where $x_1,\dots,x_n$ is the dual basis of $E$ and 
$s_1,\dots,s_n$ is the corresponding basis for the polynomial part of
the cohomology. Let $c_{ij}^k=[e_i,e_j]_k$ be the structure constants
of $\mathfrak{L}$ with respect to the basis $E$. 
We have shown that
$$
[\lambda(e_i)](e_j)=\sum_{k=1}^nc_{ij}^ke_k.
$$
Recalling that $\lambda$ is the dual of $-\beta_1 : S^* \rightarrow W^* 
\otimes S^*$, we see easily that
$$
-\beta_1(s_k) = \sum_{i,j=1}^nc_{ij}^kx_is_j.
$$
Using that $c_{ij}^k = -c_{ji}^k$, one can easily deduce the following 
equation for $\beta$:
\begin{equation}
\label{eq: exp B}
\beta(s_k)= 
 \sum_{i,j=1}^nc_{ij}^ks_ix_j + \beta_2(s_k).
\end{equation}

Thus we can also recover the Lie algebra structure of $\mathfrak{L}$
from knowledge of the $\beta_1$ component of $\beta|_{S^*}$. 
It remains to study $\beta_2$ which is the
remaining term of the Bockstein. Before doing that, we will digress to give
an application of equation~\ref{eq: exp B}.

\subsection{Comodule algebra structure} 

Let $G \in \Obj(\mathbf{LGrp})$ then the multiplication in $G$ induces
a homomorphism $\mu : \Omega_1(G) \times G \rightarrow G$ and
$$
\Delta = \mu^* : H^*(G;\mathbb{F}_p) \rightarrow
H^*(\Omega_1(G);\mathbb{F}_p) \otimes H^*(G;\mathbb{F}_p)
$$
gives $H^*(G;\mathbb{F}_p)$ the structure of a comodule algebra over
the Hopf algebra \break 
$H^*(\Omega_1(G);\mathbb{F}_p).$ (See \cite{W2}.)

In order to describe $\Delta$, it is enough to describe
it on algebra generators of $H^*(G;\mathbb{F}_p)$. To do this let us write
$$
H^*(G;\mathbb{F}_p) = \wedge^*(x_1,\dots,x_n) \otimes
\mathbb{F}_p[s_1,\dots,s_n]
$$
and $H^*(\Omega_1(G); \mathbb{F}_p) = \wedge^*(t_1,\dots,t_n) \otimes
\mathbb{F}_p[\bar{s}_1,\dots,\bar{s}_n]
$ with $\bar{s}_i = \beta(t_i)$ for
all $i$. Furthermore we have made our choices so that $j^*(s_i)=\bar{s}_i$
where $j$ is the inclusion map of $\Omega_1(G)$ into $G$ and also so that
equation~\ref{eq: exp B} holds. Recall that the
Hopf algebra structure of $\Omega_1(G)$ is one where all the generators
$\bar{s}_i$ and $t_i$ are primitive. 
Now we can state:

\begin{cor}
\label{cor: comodule}
Let $G \in \Obj(\mathbf{LGrp})$ and $c_{ij}^k$ be the structure constants of
$Log(G)$ in a suitable basis. Using the notation of the paragraph above, 
one has the following formulas which determine $\Delta$:
\begin{align*}
\Delta(x_k) =& 1 \otimes x_k \\
\Delta(s_k) =& \bar{s}_k \otimes 1 + 1 \otimes s_k + \sum_{i,j}c_{ij}^k
t_i \otimes x_j
\end{align*}
\end{cor}
\begin{proof}
Precomposing $\mu$ with the natural inclusions of $\Omega_1(G)$ and $G$
into $\Omega_1(G) \times G$, it is not hard to show
that $\Delta=\mu^*$ must take the following form on the generators:
\begin{align*}
\Delta(x_k) =& 1 \otimes x_k \\
\Delta(s_k) =& \bar{s}_k \otimes 1 + 1 \otimes s_k + \sum_{i,j}a_{ij}^k
t_i \otimes x_j.
\end{align*}
Thus it is enough to show $a_{ij}^k = c_{ij}^k$ for all $i,j,k$.
As the Bockstein $\beta$ is a natural operation, and since $\Delta$ is induced
from a homomorphism of groups, they must commute, that is $\Delta \circ
\beta = (\beta \otimes 1 + 1 \otimes \beta)\circ \Delta$. Applying 
$\beta \otimes 1 + 1 \otimes \beta$
to the formula for $\Delta(s_k)$ and using equation~\ref{eq: exp B}, we get

\begin{align}
\label{eq: icky1}
(\beta \otimes 1 + 1 \otimes \beta)(\Delta(s_k)) = 1 \otimes \beta(s_k) + 
\sum_{i,j=1}^na_{ij}^k(\bar{s_i} \otimes x_j - t_i \otimes \beta(x_j))  
\end{align}
but
\begin{align*}
\Delta(\beta(s_k)) &= \sum_{l,m=1}^nc_{lm}^k\Delta(s_lx_m) + 
\Delta(\beta_2(s_i)) \\
&= \sum_{l,m=1}^nc_{lm}^k(\bar{s}_l \otimes x_m + 1 \otimes s_lx_m
+ \sum_{u,v}a_{uv}^lt_u \otimes x_vx_m) + 1 \otimes \beta_2(s_i)
\end{align*}

Now we know the expressions in equation~\ref{eq: icky1} and the last one
are equal so equating their $H^2(\Omega_1(G);\mathbb{F}_p) \otimes
H^1(G;\mathbb{F}_p)$ components we get:
$$
\sum_{i,j=1}^na_{ij}^k\bar{s_i} \otimes x_j = \sum_{l,m=1}^nc_{lm}^k\bar{s_l}
 \otimes x_m
$$
which immediately gives $a_{ij}^k = c_{ij}^k$ for all $i,j,k$ which is what
we desired.

\end{proof}
Thus we see that the comodule algebra structure of $H^*(G;\mathbb{F}_p)$
also determines the Lie algebra corresponding to $G$.

\subsection{The class $[\mathbf{\eta}]$ and lifting uniform towers}
\label{sec: Towers}
In this section we will show that the $\beta_2$ component 
 for the Bockstein, defines a cohomology
class in a suitable cohomology group.

Fix $\mathfrak{L}=(W,br) \in \Obj(\mathbf{Lie})$. 
Recall the existence of a Koszul complex:
$$
0 \rightarrow \wedge^0(\mathfrak{L},ad) \overset{d}{\rightarrow} \dots 
\overset{d}{\rightarrow} \wedge^n(\mathfrak{L},ad) \rightarrow 0
$$
whose cohomology is $H^*(\mathfrak{L};ad)$. 
Here $\wedge^i(\mathfrak{L},ad)$ denotes
the $\mathfrak{L}$-valued alternating i-forms on $\mathfrak{L}$ and $d$ is given
by
\begin{align*}
\begin{split}
(d\omega)(u_0,\dots,u_l)=& \sum_{i<j}(-1)^{i+j}\omega 
([u_i,u_j],u_0,\dots,\Hat{u}_i,\dots,\Hat{u}_j,\dots,u_l) \\ 
& + \sum_{i=0}^l(-1)^i[u_i,\omega(u_0,\dots,\Hat{u}_i,\dots,u_l)]
\end{split}
\end{align*}
on $l$-forms $\omega$. 
If we view $\omega : \wedge^l(W) \rightarrow W$ and set
$ad : W \otimes W \rightarrow W$ to be the adjoint action map, one can
show easily that
$$
d\omega = -\omega \circ br + ad \circ (1 \wedge \omega). 
$$
Using $\btwo$ instead of $\beta_2$ in equation~\ref{eq: cohclass} and
taking duals one obtains:
$$
0 = - (-\btwo^*) \circ br + (-\bone^*) \circ (1 \wedge -\btwo^*).
$$
However, we have shown that $-\bone^*=ad$ and so we can conclude that
$-\btwo^* : \wedge^3(W) \rightarrow V$ is a closed 3-form in the Koszul
resolution for $H^*(\mathfrak{L};ad)$. In a similar fashion, when one
looks at the effect of changing the choice of the subspace $S^*$ via a
map $\mu : V^* \rightarrow \wedge^2(W^*)$ as in 
equation~\ref{eq: etaboundaryII}, one sees that the net result is
to add $d(-\mu^*)$ to $-\btwo^*$. So although the 3-form $-\btwo^*$ is
not uniquely determined from $\mathfrak{L}$, 
it is uniquely determined up to the addition of
a boundary in the Koszul resolution above 
and hence defines a unique cohomology class 
$[\eta] \in H^3(\mathfrak{L};ad)$. In fact it is easy to see that
one can get any representative of this cohomology class by
suitable choice of $\mu$. 

So now it remains to
study $[\mathbf{\eta}] \in H^3(\mathfrak{L};ad)$. It is true for $p \neq 3$ that 
$[\mathbf{\eta}]=0$ if and only if $\mathfrak{L}$ lifts to a $\mathbb{Z}/p^2\mathbb{Z}$ Lie
algebra. This requires a bit of work using results of lifting Lie
algebras as mentioned in~\cite{pak} and T. Weigel's Exp-Log
correspondence.
 We will first show that
$[\mathbf{\eta}]=0$ if and only if $G=Exp(\mathfrak{L})$
 has an $\Omega_1$-extension with
the $\Omega$EP (these terms will be defined later). This on the other hand 
is quite elementary and does not need
knowledge of the Exp-Log correspondence.

We now discuss various concepts necessary to see this fact.
First recall definitions~\ref{defn: omega1} and
\ref{defn: ppower}. 
Furthermore we define:
\begin{defn}
G is powerful if $[G,G] \subseteq G^p$.
\end{defn}
\begin{defn}
A tower of groups will be a sequence of p-central groups
$$
\{G_1,\dots,G_N\}
$$ 
(where $N=1,2,\dots,\infty$) equipped with surjective
homomorphisms $\pi_i : G_i \rightarrow G_{i-1}$ with kernel $\Omega_1(G_i)$
for $i=1,\dots, N$. ($G_0=1$ by convention.)
\end{defn}
Note the definition implies $G_1=\Omega_1(G_1)$ and hence is elementary
abelian.
Also note if we have a tower, we have central short exact sequences of groups
$$
1 \rightarrow \Omega_1(G_i) \rightarrow G_i \overset{\pi_i}{\rightarrow}
G_{i-1} \rightarrow 1
$$
and these pull back to central short exact sequences
$$
1 \rightarrow \Omega_1(G_i) \rightarrow E_i \overset{\pi_i}{\rightarrow}
\Omega_1(G_{i-1}) \rightarrow 1
$$
as considered in the beginning of this paper. Thus these have a $p$-power map
$\phi$ where $\phi(x) = \Hat{x}^p$ for $x \in \Omega_1(G_{i-1})$.
Note the $\phi$ map is automatically injective by definition of 
$\Omega_1(G_i)$. 

\begin{defn}
A uniform tower of groups is a tower of groups where all the $p$-power
maps $\phi$ of the tower are isomorphisms.
\end{defn}

Thus a uniform tower is where the $p$-power maps $\phi$ are also surjective.
Note it is well known that if $G$ is $p$-central then $G/\Omega_1(G)$
is also $p$-central so all $p$-central groups fit into a tower.
By induction, it is also easy to see that $G_n$ a group in the
$n$th stage of a uniform tower has exponent $p^n$ and is powerful.
 
In general if you have a tower, then $G_2$ fits in a short exact sequence
$$
1 \rightarrow \Omega_1(G_2) \rightarrow G_2 \rightarrow G_1 \rightarrow 1
$$
and as $G_1$ is elementary abelian, one sees immediately that $G_2$
groups in uniform towers are exactly the groups associated to bracket algebras
considered before. Now given a tower of length $N$
with end group $G_N$ one can ask when it extends to a tower of length
$N+1$, i.e., when does there
exist a $p$-central group $G_{N+1}$ such that $G_{N+1}/
\Omega_1(G_{N+1})$
is isomorphic to $G_N$. We will be primarily interested in uniform towers
and extending to uniform towers. This question has been considered
independently in \cite{W1}. In order to coincide with
the notation there we introduce the concept of $\Omega$EP used in
\cite{W1} and \cite{W2}. We will provide theorems which were proved 
independently by us
for the uniform case
and refer the reader to \cite{W1} and \cite{W2} for lifting in the
more general situation.

\begin{defn}
A $p$-central group $G$ is said to have the $\Omega$EP if there exists
a $p$-central group $\Hat{G}$ such that $\Hat{G}/\Omega_1(\Hat{G})=G$.
\end{defn}

Note that if we have such an extension 
with $\Hat{G} \overset{\pi}{\rightarrow}
G$ then $\pi^{-1}(\Omega_1(G)) \subset \Hat{G}$ is an extension of 
$\Omega_1(G)$ which has $\phi$ map injective. We call such an extension
an $\Omega_1$ extension. Let $n=\dime(\Omega_1(G))$ and 
$l=\dime(\Omega_1(\Hat{G}))$. In 
$$
H^2(\Omega_1(G);\Omega_1(\Hat{G}))=[\wedge^*(x_1,\dots,x_n) \otimes
\mathbb{F}_p[s_1,\dots,s_n]]^l
$$ 
this extension for $\pi^{-1}(\Omega_1(G))$ 
is represented by an element of the form 
$$
(s_1+\mu_1,\dots,s_n+\mu_n,\mu_{n+1},\dots,\mu_l)
$$ 
with $\mu_i \in 
\wedge^*(x_1,\dots,x_n)$. Note that $s_1+\mu_1,\dots,s_n+\mu_n$ are 
in the image of the restriction $i^* : H^*(G;\mathbb{F}_p) \rightarrow
H^*(\Omega_1(G);\mathbb{F}_p)$. Thus one has $n$ algebraically independent
elements in the second grading in the image of the restriction $i^*$.
On the other hand if we had a $p$-central group such that
$\dime(\Omega_1(G))=n$ and there are $n$ algebraically independent
elements of degree two in image $i^* : H^*(G;\mathbb{F}_p) \rightarrow
H^*(\Omega_1(G);\mathbb{F}_p)$. Then by choosing the right basis for
$$
H^2(\Omega_1(G);\mathbb{F}_p)=\wedge^*(x_1,\dots,x_n) \otimes
\mathbb{F}_p[s_1,\dots,s_n]
$$ 
one has these elements of the form
$s_1 + \mu_1,\dots,s_n+\mu_n$ with $\mu_i \in \wedge^2(x_1,\dots,x_n)$.
Thus taking an element in $H^2(G;\Omega_1(G))$ which restricts to
$(s_1+\mu_1,\dots,s_n+\mu_n) \in H^2(\Omega_1(G);\Omega_1(G))$,
we see easily that the group this element represents say $\Hat{G}$ has
$\Hat{G}/\Omega_1(\Hat{G})=G$. Thus G has the $\Omega$EP. So we conclude
\begin{lem}
\label{lem: coh-omega1}
Let $G$ be a $p$-central group with $n=\dime(\Omega_1(G))$. 
Then $G$ has the $\Omega$EP if and only if there are $n$ algebraically 
independent elements of degree 2 in image $i^* : H^*(G;\mathbb{F}_p) \rightarrow
H^*(\Omega_1(G);\mathbb{F}_p)$. 
\end{lem}

Now note in a tower $\{G_1,\dots,G_N\}$, $G_i$ has the $\Omega$EP for
$1 \leq i < N$ by definition. This lets us prove:
\begin{thm}
\label{thm: towercoh}
Let $\{G_1,\dots,G_N\}$ be a uniform tower with $\dime(G_1)=n$. Then
for $1 \leq i < N$ we have
$$
H^*(G_i;\mathbb{F}_p)=\wedge^*(x_1,\dots,x_n) \otimes \mathbb{F}_p[s_1,\dots,s_n]
$$
where $\degr(x_i)=1$ (and for $i>1$, $x_i$ is induced from $G_{i-1}$) and 
$\degr(s_i)=2$ (and for $i>1$, $s_i$ ``comes'' from $\Omega_1(G_i)$),
and $H^*(G_N;\mathbb{F}_p)$ is given by the same formula if and only if
the tower extends to a uniform tower $\{G_1,\dots,G_N,G_{N+1}\}$.
\end{thm}
\begin{proof}
By the comments in the preceding paragraph and lemma~\ref{lem: coh-omega1}
one obtains the first part of this theorem by induction on $i$ where 
$1 \leq i <N$. More explicitly note that one has the result for $G_1$. So fix
$i > 1$ and assume we have shown the first part for all $j < i$. Then we
have the central short exact sequence:
$$
1 \rightarrow \Omega_1(G_i) \rightarrow G_i \overset{\pi_i}{\rightarrow}
G_{i-1} \rightarrow 1.
$$
So the $E_2^{*,*}$ term of the L.H.S.-spectral sequence of this 
extension is given by 
$$
\wedge^*(t_1,\dots,t_n) \otimes \mathbb{F}_p[s_1,\dots,s_n]
\otimes \wedge^*(x_1,\dots,x_n) \otimes \mathbb{F}_p[y_1,\dots,y_n]
$$
where the first two terms are the cohomology of $\Omega_1(G_i)$ and the
last two terms are the cohomology of $G_{i-1}$ by induction. One has
as usual $d_2(s_i)=d_2(x_i)=d_2(y_i)=0$ and $d_2(t_i)$ are the elements
representing the extension. By the proof of lemma~\ref{lem: coh-omega1} 
and the fact $G_i$ is an $\Omega_1$ extension of $G_{i-1}$ one sees (after
changing basis for $y$'s) that
these elements are of the form $d_2(t_i)=y_i + \mu_i$ for all $i$ where
$\mu_i \in \wedge^*(x_1,\dots,x_n)$. Thus as in the proof of 
theorem~\ref{thm: Fp cohomology}, 
one gets $E_3^{*,*}=\mathbb{F}_p[s_1,\dots,s_n] \otimes
\wedge^*(x_1,\dots,x_n)$ with $d_3(x_i)=0$. Now $E_{\infty}^{0,*}$
is the image of $i^* : H^*(G_i;\mathbb{F}_p) \rightarrow H^*(\Omega_1(G_{i});\mathbb{F}_p)$ and
by lemma~\ref{lem: coh-omega1} this contains $n$ algebraically independent
elements in degree 2. Since $E_2^{0,2}$ has dim $n$ we see by necessity
that $d_3=0$ on $E_2^{0,2}$, which implies $d_3=0$. Then
as in theorem~\ref{thm: Fp cohomology} one sees that $H^*(G_i;\mathbb{F}_p)$
is as stated. So we are done with the first part. If the tower extends to
a uniform tower of length one more, then one can perform the same argument
on $G_N$ to conclude its cohomology is as stated. On the other hand if its
cohomology is as stated, image $i^* : H^*(G_N;\mathbb{F}_p) \rightarrow
H^*(\Omega_1(G_N);\mathbb{F}_p)$ has $n$ algebraically independent elements
in degree 2 and so $G_N$ has the $\Omega$EP. Using these elements as in
the proof of lemma~\ref{lem: coh-omega1}, one constructs an extension to
a uniform tower of length $N+1$.
\end{proof}

\textbf{Examples:}
\begin{defn}
Fix $p$ an odd prime number, then
$\Gamma_{n,k}(p)$ is defined as the kernel of the reduction homomorphism
$$
GL_n(\mathbb{Z}/p^{k+1}\mathbb{Z}) \rightarrow
GL_n(\mathbb{Z}/p\mathbb{Z})
$$
for all $n,k \geq 1$. Similarly, $\Hat{\Gamma}_{n,k}(p)$ is the group obtained
by replacing $GL_n$ with $SL_n$ in the definition above.
\end{defn}
One shows by induction on $k$ that the $\Gamma_{n,k}(p)$ 
 are $p$-groups and it is then easy to show that for
 fixed $n$ they fit together to give an infinite uniform tower:
$$
\rightarrow \Gamma_{n,k}(p) \rightarrow \Gamma_{n,k-1}(p) \rightarrow
\dots \rightarrow \Gamma_{n,1}(p) \rightarrow 1.
$$
A similar statement holds for the $\Hat{\Gamma}_{n,k}(p)$.

Theorem~\ref{thm: towercoh} has the following immediate
corollary:
\begin{cor}
For $n, k \geq 1$ and $p$ an odd prime we have:
$$
H^*(\Gamma_{n,k}(p) ; \mathbb{F}_p) \cong 
\wedge^*(x_{1},\dots,x_{n^2}) \otimes
\mathbb{F}_p[s_{1},\dots,s_{n^2}]
$$
and
$$
H^*(\Hat{\Gamma}_{n,k}(p) ; \mathbb{F}_p) \cong
\wedge^*(x_1,\dots,x_{n^2-1}) \otimes
\mathbb{F}_p[s_1,\dots,s_{n^2-1}]
$$
where $\degr(x_i)=1$ and $\degr(s_i)=2$ for all $i$.
\end{cor}

Recall that if 
$G \in \Obj(BGrp)$ then $G=G_2$ of a uniform tower $\{G_1,G_2\}$. 
This tower extends to a uniform tower of length 3 if and only if
$Log(G)$ is a Lie algebra. This follows form theorem~\ref{thm: towercoh}
and theorem~\ref{thm: Fp cohomology}. Next, we will show that
this uniform tower $\{G_1,G_2\}$ extends to a tower of length 4 if and only if
$[\mathbf{\eta}]=0$ where $[\mathbf{\eta}] \in H^3(\mathfrak{L};ad)$ is 
the cohomology class defined earlier.

Now suppose we have a uniform tower $\{G_1,G_2,G_3'\}$ so $G_2 \in \Obj(LGrp)$
and has $\mathbb{F}_p$-cohomology given by theorem~\ref{thm: Fp cohomology}.
We will now show that there is a uniform tower $\{G_1,G_2,G_3,G_4\}$
if and only if $[\mathbf{\eta}]=0$ where we say $G_3$ instead of $G_3'$ since we
might have to choose a different $\Omega_1$ extension of $G_2$ to do this.
Note the tower $\{G_1,G_2\}$ extends to one of length 4 as desired if and
only if $G_2$ has a uniform $\Omega_1$ extension $G_3$ which itself has
the $\Omega$EP. By theorem~\ref{thm: towercoh}, this happens if and only
if $G_3$ has $\mathbb{F}_p$-cohomology $\wedge^*(x_1,\dots,x_n) \otimes
\mathbb{F}_p[s_1,\dots,s_n]$ where $\degr(x_i)=1,\degr(s_i)=2$ and $n=\dime(G_1)$.

So let us study the cohomology of an $\Omega_1$ extension $G_3$ of $G_2$.
So $H^*(G_2;\mathbb{F}_p)=\wedge^*(x_1,\dots,x_n) \otimes \mathbb{F}_p[s_1,\dots,
s_n]$ where $\degr(s_i)=2,\degr(x_i)=1$ as usual. Then $G_3$ is represented
in $H^2(G_2;\Omega_1(G_3))=[H^2(G_2;\mathbb{F}_p)]^n$ by $(s_1 + \mu_1,\dots,
s_n + \mu_n)$ where $\mu_i \in \wedge^2(x_1,\dots,x_n)$ under a suitable
choice of basis of $\Omega_1(G_3)$. However now note that due to the
ambiguous nature of $s_i$ we can shift them by elements of 
$\wedge^2(x_1,\dots,x_n)$ so without loss of generality the extension element
of $G_3$ is $(s_1,\dots,s_n)$. Note this implies a particular choice of
$s$-basis. Let $\mathbf{\eta}$ be the Bockstein term for this choice of
$s$-basis. Then as we saw in the first sections, when we look at the 
L.H.S.-spectral
sequence for
$$
1 \rightarrow \Omega_1(G_3) \rightarrow G_3 \overset{\pi_3}{\rightarrow}
G_2 \rightarrow 1
$$
Then 
$$
E_2^{*,*}=\wedge^*(t_1,\dots,t_n) \otimes \mathbb{F}_p[\beta(t_1),\dots,
\beta(t_n)]
\otimes \wedge^*(x_1,\dots,x_n) \otimes \mathbb{F}_p[s_1,\dots,s_n]
$$ 
where the first two terms are the cohomology of $\Omega_1(G_3)$  
and the last
two are that of $G_2$. Again $d_2(\beta(t_i))=d_2(x_i)=d_2(s_i)=0$ and
$d_2(t_i) = s_i$. 
Thus,
$$
E_3^{*,*}=\mathbb{F}_p[\beta(t_1),\dots,\beta(t_n)]
\otimes \wedge^*(x_1,\dots,x_n)
$$ 
with $d_3(x_i)=0$ and
$d_3(\beta(t_i))=\{ \beta(d_2(t_i)) \}=\{ \beta(s_i) \}$ in $E_3$. 
So by what we
said before $G_3$ has the $\Omega$EP if and only if $\{ \beta(s_i) \}=0$ for all
$i$. So let us calculate this explicitly. There are elements
$\xi_{i,l} \in \wedge^1(x_1,\dots,x_n)$ such that: 
\begin{align*}
\beta(s_i) = \sum_{l=1}^n\xi_{i,l}s_l + \beta_2(s_i)
	 = \sum_{l=1}^n\xi_{i,l}(0) + \beta_2(s_i) \\ 
\end{align*}
where in the last step we used that $s_i=0$ in $E_3$ due
to the fact $d_2(t_i) = s_i$. 
Thus we get $d_3(\beta(t_i))=\beta_2(s_i)$ for all $i$. So that $G_3$ defined by the
extension element $(s_1,\dots,s_n)$ has the $\Omega$EP if and only if
$\beta_2=0$. As we have argued that the $\Omega_1$ extensions of $G_2$
are exactly those defined by such extension elements, we see that there
is an $\Omega_1$ extension $G_3$ of $G_2$ such that $G_3$ has the $\Omega$EP
if and only if there is a choice of $s_1,\dots,s_n$ via changing up
to adding elements in $\wedge^2(x_1,\dots,x_n)$ so that the corresponding
$\beta_2=0$, i.e., if and only if
$[\mathbf{\eta}]=0$. Thus we have shown:
\begin{thm}
\label{thm: eta lift}
Let $\{G_1,G_2\}$ be a uniform tower. Then it extends to a uniform tower
of length 3 if and only if $\mathfrak{L}=Log(G_2)$ is a Lie algebra. In this case
it is possible to extend $\{G_1,G_2\}$ to a uniform tower of length 4
if and only if $[\mathbf{\eta}] \in H^3(\mathfrak{L};ad)$ is zero.
\end{thm}

\subsection{Formulas for $B_2^*$}
In this section we will derive explicit formulas for $B_2^*$ of the
Bockstein Spectral Sequence of a group in $\mathbf{LGrp}$. These
formulas will express $B_2^*$ in terms of various Lie algebra cohomologies
of the Lie algebra corresponding to the group.

Let $G \in \Obj(\mathbf{LGrp})$ and $Log(G)=\mathfrak{L}$. Let
$n=\dime(\mathfrak{L})$ and assume $[\mathbf{\eta}]=0$. Then we have that
$$
H^*(G;\mathbb{F}_p) = \wedge^*(x_1,\dots,x_n) \otimes \mathbb{F}_p[s_1,\dots,s_n].
$$
The Bockstein is given by 
$$
\beta(x_i)=-\sum_{l < m}c_{lm}^ix_lx_m
$$
and
$$
\beta(s_i)=\sum_{l,m=1}^nc_{lm}^is_lx_m
$$
where $c_{lm}^i$ are the 
structure constants of $\mathfrak{L}$ in a basis $\{e_1,\dots,e_n\}$
implicit in all this. (Note the other Steenrod
operations are axiomatically determined so indeed we have determined
$H^*(G;\mathbb{F}_p)$ as a Steenrod module. So we see the Steenrod module
$H^*(G;\mathbb{F}_p)$ determines the group in the class of objects $\mathbf{LGrp}$.)
Thus we can view $H^*(G;\mathbb{F}_p)$ as being bigraded in the same
way as $E_{\infty}^{*,*}$ of the LHS-spectral sequence which we used
to compute it.
We see that when $[\mathbf{\eta}]=0$, the Bockstein $\beta$
preserves the polynomial degree in this bigrading.
Thus the complex $(B_1^*,\beta=\beta_1)$
of the Bockstein spectral sequence splits up as a direct sum of
differential complexes, one for each polynomial degree. 
More precisely, we have 
$$
\mathbf{F}_p[s_1,\dots,s_n] = \oplus_{k=0}^{\infty}S^k
$$
where $S^k$ stands for the vector space of homogeneous polynomials of
degree $k$. Thus we have
$$
H^*(G;\mathbf{F}_p)= \oplus_{k=0}^{\infty}(\wedge^*(x_1,\dots,x_n)
\otimes
S^k)
$$
where each $\wedge^*(x_1,\dots,x_n) \otimes S^k$ becomes a
differential
complex under $\beta$. Thus $B_2^*$ is
given as the direct sum of the cohomology of each of these complexes.
We can identify these complexes with the Koszul complexes used to
calculate $H^*(\mathfrak{L};S^k)$ where $S^k$ is the vector space of
symmetric multilinear k-forms on $\mathfrak{L}$ given the Lie algebra action:
$$
(u.f)(u_1,\dots,u_k)=\sum_{i=1}^kf(u_1,\dots,[u_i,u],\dots,u_k)
$$
for all $f \in S^k$ and $u \in \mathfrak{L}$. 
To see this note that the direct sum of these Koszul complexes has
a natural identification with
$$
\oplus_{k=0}^{\infty}(\wedge^*(x_1,\dots,x_n) \otimes S^k)
$$
and hence $H^*(G;\mathbf{F}_p)$. In \cite{pak}, it is
shown that the sum of the Koszul differentials is a derivation with
respect to the algebra structure given, i.e., that of
$H^*(G;\mathbf{F}_p)$.
So to show that these are the same as the Bockstein differential it
is enough to check that they agree on the generators $x_1, \dots,
x_n, s_1, \dots,s_n$ of $H^*(G;\mathbf{F}_p)$. This is easy and will be left to
the reader (for details see \cite{pak}). Thus once, grading is taken into
consideration one gets the equation:
\begin{equation}
B_2^* = \oplus_{k=0}^{\infty}H^{*-2k}(\mathfrak{L},S^k)
\end{equation}
which holds when $[\mathbf{\eta}]=0$.
Thus we see that knowing the second term of the Bockstein
spectral sequence is equivalent to knowing all the Lie algebra cohomology
groups $H^*(\mathfrak{L};S^k)$. For clarity we state the results of
this section in the following theorem.
\begin{thm}
\label{thm: B2}
Let $G \in \Obj(\mathbf{LGrp})$ with $Log(G)=\mathfrak{L}$ and with $[\mathbf{\eta}]=0$.
Then $B_2^*$ of the Bockstein spectral sequence for $G$ is given by
$$
B_2^*=\oplus_{k=0}^{\infty}H^{*-2k}(\mathfrak{L};S^k).
$$
\end{thm} 
We will calculate $B_2^*$ (at least partially) in some cases 
later on in this paper.

\section{The General Case}
\subsection{Introduction}
Our results so far concern groups $G$ which
fit in as $G_2$ in some uniform tower 
$$
\dots \rightarrow G_n \rightarrow \dots \rightarrow G_2 \rightarrow
G_1 \rightarrow 1
$$
We now study the higher groups $G_n$ in such a tower. The
$\mathbb{F}_p$-cohomology
of such groups is determined by theorem~\ref{thm: towercoh}.
We will now extend our formulas
for the Bockstein to the higher groups $G_n$ , $n \geq 3$ in 
such uniform towers.

To do this we will need to use T. Weigel's Exp-Log correspondence.
Let $\mathbb{Z}_p$ denote the $p$-adic integers. 
\begin{defn}
A $\mathbb{Z}_p$-Lie algebra is a $\mathbb{Z}_p$-module $\mathfrak{L}$ equipped
with a bilinear, alternating product $[\cdot,\cdot] : \mathfrak{L}
\times \mathfrak{L} \rightarrow \mathfrak{L}$ which satisfies the 
Jacobi Identity.  
\end{defn}
\begin{defn}
A $\mathbb{Z}_p$-Lie algebra $\mathfrak{L}$ is powerful if
$[\mathfrak{L},\mathfrak{L}] \subseteq p\cdot \mathfrak{L}$.
\end{defn}
\begin{defn}
$\Omega_1(\mathfrak{L})= \{ v \in \mathfrak{L} : p\cdot v=0 \}$
\end{defn}
\begin{defn}
A $\mathbb{Z}_p$-Lie algebra $\mathfrak{L}$ is p-central if
$[\mathfrak{L},\Omega_1(\mathfrak{L})]=0$.
\end{defn}
\begin{defn}
$\mathbf{CLie}$ is the category of powerful, p-central finite
$\mathbb{Z}_p$-Lie algebras, with morphisms the maps of Lie algebras.
\end{defn}
\begin{defn}
$\mathbf{CGrp}$ is the category of powerful, p-central p-groups
with morphisms the homomorphisms of groups.
(see section~\ref{sec: Towers} for the relevant definitions.)
\end{defn}
\begin{defn}
In analogy to the case with groups we say that $\mathfrak{L} \in
\Obj(\mathbf{CLie})$ has the $\Omega$EP if there exists
$\bar{\mathfrak{L}} \in \Obj(\mathbf{CLie})$ such that
$\bar{\mathfrak{L}}/\Omega_1(\bar{\mathfrak{L}}) \cong \mathfrak{L}$.
\end{defn}
The following theorem was proved in \cite{W3} and stated in
\cite{W1}. 
It generalizes the correspondence
between p-adic Lie algebras and uniform pro-p groups (see \cite{DS}) 
and is based on a careful study of the Baker-Campbell-Hausdorff
identity. ($e^A \cdot e^B = e^C$ where $C$ lies in the Lie algebra
generated by $A$ and $B$.)
\begin{thm}{\cite{W3}: Exp-Log correspondence.}
\label{thm: exp-log}
Let $p \geq 5$. There exist functors 
$\mathbf{Exp} : \mathbf{CLie} \rightarrow \mathbf{CGrp}$ and
$\mathbf{Log} : \mathbf{CGrp} \rightarrow \mathbf{CLie}$ which
give a natural equivalence of categories. Furthermore if
$\mathfrak{L} \in \Obj(\mathbf{CLie})$ then $\mathfrak{L}$ has
the $\Omega$EP if and only if $\mathbf{Exp}(\mathfrak{L})$ does.
\end{thm}
Note that this theorem does not cover the case $p=3$.

\begin{defn}
$\mathbf{CLie}_k$ is the full subcategory of $\mathbf{CLie}$
whose objects are the powerful, p-central Lie algebras which are
free $\mathbb{Z}/p^k\mathbb{Z}$-modules.
\end{defn}
\begin{defn}
$\mathbf{CGrp}_k$ is the
full subcategory of $\mathbf{CGrp}$ consisting of groups which
fit into the kth stage of some uniform tower.
\end{defn} 
The ``exp-log'' equivalence above,
can easily be shown to restrict to a natural equivalence of
$\mathbf{CGrp}_k$ with $\mathbf{CLie}_k$ for all $k \geq 1$.
\begin{defn}
$\mathbf{Brak}_{k}$ is the category whose objects consist of
bracket algebras over $\mathbb{Z}/p^k\mathbb{Z}$. That is
free $\mathbb{Z}/p^k\mathbb{Z}$-modules of finite rank $B$, equipped
with a bilinear, alternating map $[\cdot,\cdot] : B \times B
\rightarrow B$. The morphisms are maps of the underlying modules which
preserve the bracket structure.
\end{defn}
In \cite{pak},
we construct a shifting function $\mathbf{S} : \mathbf{CLie}_k
\rightarrow \mathbf{Brak}_{k-1}$ which is a bijection on the object
level and induces certain identifications on the morphisms. We also
show that for $\mathfrak{L} \in \Obj(\mathbf{CLie}_k)$, then 
$\mathfrak{L}$ has the $\Omega$EP if and only if
$\mathbf{S}(\mathfrak{L})$ is a Lie algebra. Let $Log :
\mathbf{CGrp}_k \rightarrow \mathbf{Brak}_{k-1}$ be the composition of
the functors $\mathbf{Log}$ and $\mathbf{S}$. Note this extends the
definition of $Log: \mathbf{CGrp}_2 \rightarrow \mathbf{Brak}_1$
studied earlier. The facts mentioned
above, together with theorem~\ref{thm: exp-log} allows one to prove
the following theorem on uniform towers:

\begin{thm}
\label{thm: Tower extension theorem}
Fix $p \geq 5$. Let
$$
G_k \rightarrow \dots \rightarrow G_2 \rightarrow G_1 \rightarrow 1
$$
be a uniform tower of groups, of length $k$. 
Then this tower extends to a uniform
tower of length $k+1$ if and only if $Log(G_k)$ is a Lie algebra.
\end{thm}
\begin{proof}
The tower extends to a uniform tower of length $k+1$ if and only if
$G_k$ has the $\Omega$EP. By the comments made before, this happens
if and only if $Log(G_k)$ is a Lie algebra.
\end{proof}
Note this extends the $k=2$ case proven in theorem~\ref{thm: Fp cohomology}. 
Let $k \geq 2$.
Since $Log: \mathbf{CGrp}_k \rightarrow \mathbf{Brak}_{k-1}$ is
a bijection on the object level, we will refer to $Exp(B)$ as the
unique group $G \in \Obj(\mathbf{CGrp}_k)$ such that $Log(G)=B$.
If one has a morphism $\psi \in \Mor_{\mathbf{Brak}_{k-1}}(B_1,B_2)$,
then one gets a map of groups $\Psi : Exp(B_1) \rightarrow Exp(B_2)$
which is unique up to multiplication by a map $\mathbf{\eta} : Exp(B_1)
\rightarrow \Omega_1(Exp(B_2))$. (see \cite{pak}.)
Although this map is not unique, it induces a unique group
homomorphism on the lower levels of the tower. ( $Exp(B_i)$ fits into
the kth stage of some uniform tower.) 

\begin{defn}
For this section, we fix $p \geq 5$, a prime.
Let $\mathbf{LGrp}_k$ denote the full subcategory of $\mathbf{CGrp}_k$
consisting of the objects $G \in \Obj(\mathbf{CGrp}_k)$ which have
$Log(G)$ a Lie algebra. Let $\mathfrak{L}=(W,br)$ be the reduction
of $Log(G)$ to a $\mathbb{F}_p$-Lie algebra.
\end{defn}
Then by theorems~\ref{thm: towercoh} and \ref{thm: Tower extension
theorem}, one has that
$$
H^*(G ; \mathbb{F}_p) \cong \wedge^*(x_1,\dots,x_n) \otimes
\mathbb{F}_p[s_1,\dots,s_n]
$$
where $n = \dime(\mathfrak{L})$. Recall the exterior algebra comes 
isomorphically from the lower level in the tower so can be identified
with the exterior algebra at the first stage of the tower. 
Thus by naturality of the Bockstein, we see that
$$
\beta(x_i) = -br^*(x_i)
$$
which is the same formula as in the second stage of the tower.
The
polynomial algebra part maps isomorphically to that of $\Omega_1(G)$
which is isomorphic to the first level of the tower as a group via
an iterate of the pth power map. We will use this isomorphism to
choose basis from now on without mention - thus maps on the group
level will induce the same maps on the $x_i$ as they do on the $s_i$.
(Up to addition of exterior elements.)
As before, we reason that
$$
\beta(\mathbf{s}) = \xi \mathbf{s} + \mathbf{\eta}
$$
where $\mathbf{\eta}$ is a column vector of elements in
$\wedge^3(x_1,\dots,x_n)$,
$\xi$ is a $n \times n$ matrix of elements in
$\wedge^1(x_1,\dots,x_n)$
and $\mathbf{s}$ is the column vector of the polynomial generators.
As before, because $\beta \circ \beta = 0$, one can show $\xi$
determines a map of Lie algebras $\lambda : \mathfrak{L} \rightarrow 
\mathfrak{gl}(\mathfrak{L})$.

Let $G_1, G_2 \in \Obj(\mathbf{LGrp}_k)$ and
$\psi \in \Mor(Log(G_1),Log(G_2))$. Let $\mathfrak{L}_i$ be the 
$\mathbb{F}_p$-Lie algebras obtained by reducing $Log(G_i)$ for
$i=1,2$ and let $\bar{\psi} : \mathfrak{L}_1 \rightarrow
\mathfrak{L}_2$ be the map of Lie algebras induced by $\psi$. 
Then $Exp(\psi) : G_1 \rightarrow G_2$ induces a map on the
cohomology level which is well defined and induced by pullback
via $\bar{\psi}$ on the exterior generators and which is given
by the same formula on the polynomial generators except it is ambiguous
up to addition of exterior elements. As before this enables us to
prove naturality of the $\lambda$ maps, that is the various
$\lambda$ maps fit together to form a strong natural self-representation
$\lambda$ on $\mathbf{Lie}_{k-1}(p)$. 

By theorem~\ref{thm: ad or 0 II}, to show that $\lambda = ad$ it is enough to show that it is
nontrivial
on $\mathfrak{sl}_2$ and $\mathfrak{so}_3$. This follows exactly as it
did before by considering the Bockstein spectral sequence and noting
that it must converge to zero in positive gradings. Thus if we
fix $G \in \Obj(\mathbf{Lie}_k)$ and a basis $\{ e_1,\dots,e_n \}$ for 
$\mathfrak{L}$, and let $c_{ij}^k$
be the structure constants for $\mathfrak{L}$ with respect to this
basis, 
then we obtain the following formula for the Bockstein
as before:
\begin{align*}
\begin{split}
\beta(x_k) =& -\sum_{i<j} c_{ij}^kx_ix_j \\
\beta(s_k) =& \sum_{i,j=1}^n c_{ij}^ks_ix_j + \eta_k
\end{split}
\end{align*}
Again, one shows $\mathbf{\eta}$ determines a well-defined cohomology class
$[\mathbf{\eta}] \in H^3(\mathfrak{L};ad)$ which vanishes if and only if 
the uniform tower that $G$ is in can be extended two steps above $G$.
We summarize in the following theorem:
\begin{thm}
\label{thm: main}
Fix $p \geq 5$. Let
$$
G_{k+1} \rightarrow G_k \rightarrow \dots \rightarrow G_2 \rightarrow
G_1 \rightarrow 1
$$
be a uniform tower with $k \geq 2$. Let $\mathfrak{L}=Log(G_2)$
and let $c_{ij}^k$ be the structure constants of $\mathfrak{L}$ with
respect to some basis. Then for suitable
choices of degree 1 elements $x_1,\dots,x_n$ and degree 2
elements $s_1,\dots,s_n$, one has
$$
H^*(G_k ; \mathbb{F}_p) \cong \wedge^*(x_1,\dots,x_n) \otimes
\mathbb{F}_p[s_1,\dots,s_n]
$$
with 
\begin{align*}
\begin{split}
\beta(x_t) =& -\sum_{i < j}^n c_{ij}^tx_ix_j \\
\beta(s_t) =& \sum_{i,j=1}^n c_{ij}^ts_ix_j + \eta_t
\end{split}
\end{align*}
for all $t = 1,\dots n$. Furthermore $\mathbf{\eta}$ defines a 
cohomology class $[\mathbf{\eta}] \in H^3(\mathfrak{L} ; ad)$ which vanishes if
and only if there exists a uniform tower
$$
G_{k+2}' \rightarrow G_{k+1}' \rightarrow G_k \rightarrow
\dots \rightarrow G_2 \rightarrow G_1 \rightarrow 1
$$ or in other words if and only if $Log(G_k)$ (which is a Lie algebra over
$\mathbb{Z}/p^{k-1}\mathbb{Z}$), has a lift to a Lie algebra
over $\mathbb{Z}/p^k\mathbb{Z}$. In this case, one can drop the $\mathbf{\eta}$
terms in the formula for the Bockstein.
\end{thm}

Note, for an example of a Lie algebra over $\mathbb{Z}/p^{k-1}\mathbb{Z}$ 
which does
not lift to a Lie algebra over $\mathbb{Z}/p^k\mathbb{Z}$ see~\cite{Br}.

Note that the obstruction classes $[\mathbf{\eta}]$ are functionally equivalent
to the obstruction classes $J_{[\cdot,\cdot]}$ discussed in
\cite{Br}. Also note that if 
$$
H^*(G; \mathbb{F}_p)
\cong \wedge^*(x_1,\dots,x_n) \otimes \mathbb{F}_p[s_1,\dots,s_n]
$$
then all the other (non-Bockstein) Steenrod operations are
axiomatically determined. So the only difference in the
$\mathbb{F}_p$-cohomologies as Steenrod algebras, of the different
groups in a uniform tower (except at an end) is in the classes $[\mathbf{\eta}]$.
Thus we conclude the following corollary:
\begin{cor}
Fix $p \geq 5$. Let
$$
G_{k+2} \rightarrow G_{k+1} \rightarrow G_k \rightarrow \dots
\rightarrow G_2 \rightarrow G_1 \rightarrow 1
$$
be a uniform tower with $k \geq 2$. Then $H^*(G_i ; \mathbb{F}_p)$
are isomorphic Steenrod modules for $2 \leq i \leq k$.
\end{cor}

Finally, one should mention that theorem~\ref{thm: main} extends to  
 the category of powerful, $p$-central $p$-groups
with the $\Omega$EP. As mentioned before, there is a natural equivalence
$\mathbf{Log}$ from this category to the category
of powerful, $p$-central, finite $\mathbb{Z}_p$-Lie algebras with the 
$\Omega$EP. The shifting functor $\mathbf{S}$ mentioned before (see
also \cite{pak}) maps this category into the category of
finite $\mathbb{Z}_p$-Lie algebras. (Again, the Jacobi identity holds
exactly because of the $\Omega$EP as mentioned earlier in this paper.)
Unfortunately, the shifting functor does lose a bit of information on the
object level which corresponds to an elementary abelian summand
in the group/Lie algebra. 

Let $Log=\mathbf{S} \circ \mathbf{Log}$,
then for a powerful, $p$-central, $p$-group $G$ with the $\Omega$EP, 
one has a corresponding Lie algebra $\mathfrak{L}=Log(G)$ 
and this induces a Lie algebra over
$\mathbb{F}_p$, $\bar{\mathfrak{L}}=\mathfrak{L}/p\mathfrak{L}$.
Due to the loss of the elementary abelian summand, one can have
$\dime(\bar{\mathfrak{L}}) < \dime(\Omega_1(G))$. This is unavoidable,
as the Bockstein on the generators of the elementary abelian summand
are of a different nature (Bockstein of degree 1 generators is no longer
nilpotent for example) but easy to describe as we know what the
Bockstein is on elementary abelian $p$-groups. Taking this into account,
one can obtain the following extension of theorem~\ref{thm: main}.

\begin{thm}
\label{thm: last}
Fix $p \geq 5$. Let $G$ be a powerful, $p$-central $p$-group with
the $\Omega$EP. Let $\mathfrak{L}=Log(G)$, $\bar{\mathfrak{L}}=
\mathfrak{L}/p\mathfrak{L}$ and let $c_{ij}^k$ be the structure constants of 
$\bar{\mathfrak{L}}$ with respect to some basis. Suppose 
$\dime(\Omega_1(G))=n$ and $\dime(\bar{\mathfrak{L}})=k$. Then for suitable
choices of degree 1 elements $x_1,\dots,x_n$ and degree 2
elements $s_1,\dots,s_n$, one has
$$
H^*(G ; \mathbb{F}_p) \cong \wedge^*(x_1,\dots,x_n) \otimes
\mathbb{F}_p[s_1,\dots,s_n]
$$
with 
\begin{align*}
\begin{split}
\beta(x_t) =& -\sum_{i < j}^k c_{ij}^tx_ix_j \\
\beta(s_t) =& \sum_{i,j=1}^k c_{ij}^ts_ix_j + \eta_t \text{ for }
t = 1, \dots k. \\
\beta(x_t) =& s_t \text{ and }
\beta(s_t) = 0 \text{ for } t > k.
\end{split}
\end{align*}
Furthermore the $\mathbf{\eta}_t$ define a 
cohomology class $[\mathbf{\eta}] \in H^3(\bar{\mathfrak{L}} ; ad)$ which vanishes if and only if the Lie algebra $\mathfrak{L}$ has the $\Omega$EP. When
this is the case, one can drop the $\eta$ terms in the formula above.
\end{thm}

Here the definition of $\Omega$EP has been extended to finite 
$\mathbb{Z}_p$-Lie algebras in general by saying that $\mathfrak{L}$ has the 
$\Omega$EP if there is a $\mathbb{Z}_p$-Lie algebra $\bar{\mathfrak{L}}$ such 
that $\bar{\mathfrak{L}}/\Omega_1(\bar{\mathfrak{L}})=\mathfrak{L}$.
  
\subsection{Calculating $B_2^*$ in various cases}
Recall that there is a natural equivalence between the
categories $\mathbf{Lie}$ and $\mathbf{LGrp}$. Under this equivalence,
 Abelian Lie algebras correspond to abelian
groups so are not so interesting. Let us work out $B_2^*$ for the
nonabelian Lie algebra $\mathfrak{S}$ which has basis $\{x ,y\}$ and
bracket given by the relation $[y,x]=y$. Let the corresponding group
in $\mathbf{LGrp}$ be denoted $G(\mathfrak{S})$. This is a group of
order $p^4$ and exponent $p^2$. Then we have the structure
constants $c_{y,x}^y=1$ and $c_{y,x}^x=0$. Thus using
theorem~\ref{thm: main}, we see that
$$
H^*(G(\mathfrak{S}) ; \mathbb{F}_p) \cong \wedge^*(x,y) \otimes 
\mathbb{F}_p[X,Y]
$$
with
\begin{align*}
\begin{split}
\beta(x)= 0 \text{ , }& \beta(y)= xy \\
\beta(X)= 0 \text{ , }& \beta(Y)= Yx - Xy
\end{split}
\end{align*}
where we used that the class $[\mathbf{\eta}]$ vanishes as $\mathfrak{S}$
lifts to a Lie algebra over the $p$-adic integers.
Let $A$ be the subalgebra of $H^*(G(\mathfrak{S}) ; \mathbb{F}_p)$
generated by $x,yY^{p-1}, X,Y^p$. This is isomorphic as graded
algebras, to the graded algebra
$$
\wedge^*(x,z) \otimes \mathbb{F}_p[X,Z]
$$
where $\degr(x)=1,\degr(z)=2p-1,\degr(X)=2,\degr(Z)=2p$.
It is easy to check that $\beta$ vanishes on $A$. Thus there is a well
defined map of graded-algebras from $A$ to $B_2^*$. A direct
calculation shows that this is an isomorphism. By comparisons (restricting
to suitable subgroups of $G(\mathfrak{S})$),
it is easy to show then that we can choose $z,Z$ so that $\beta_2$
is given by $\beta_2(x)=X,\beta_2(z)=Z$. Thus $B_3^*=0$ for $* > 0$.
Thus we can conclude that 
$$
\expo(\bar{H^*}(G(\mathfrak{S}) ; \mathbb{Z}))=p^2.
$$
Now let us consider the Lie algebra $\mathfrak{sl}_2$ with basis
$\{h, x_+, x_-\}$ with bracket given by $[h,x_+]=2x_+$,$[h,x_-]=-2x_-$
and $[x_+,x_-]=h$. Let $G(\mathfrak{sl}_2)=\Hat{\Gamma}_{2,2}(p)$
be the corresponding group. Again we note this Lie algebra lifts to the
$p$-adics so by theorem~\ref{thm: main} we conclude:
$$
H^*(G(\mathfrak{sl}_2) ;\mathbb{F}_p) \cong \wedge^*(h,x_+,x_-)
\otimes \mathbb{F}_p[H,X_+,X_-]
$$
with Bockstein given by
\begin{align*}
\begin{split}
\beta(h)= -x_+x_- \text{ , }& \beta(H)= X_+x_- - X_-x_+ \\
\beta(x_+)= -2hx_+ \text{ , }& \beta(X_+)= 2(Hx_+ - X_+h) \\
\beta(x_-)= 2hx_- \text{ , }& \beta(X_-)= -2(Hx_- -X_-h).
\end{split}
\end{align*}
Using this explicit formula and viewing $B_2^*$ as the sum of Lie
algebra cohomologies as was shown before, one computes $H^*(\mathfrak{sl}_2
;\mathbb{F}_p)$ is zero except in dimensions 0 and 3 where it is one 
dimensional. Also $H^*(\mathfrak{sl}_2 ; S^1)=0$ and
$H^0(\mathfrak{sl}_2 ; S^2)$ is one dimensional generated by the
Killing form of $\mathfrak{sl}_2$. This implies the following facts:
$\dime(B_2^1)=\dime(B_2^2)=0$, $\dime(B_2^3)=1$ generated by $hx_+x_-$
and $\dime(B_2^4)=1$ generated by the Killing form $8(H^2 + X_+X_-)$.
By comparisons, $\beta_2(hx_+x_-)=0$ so that $B_2^3=B_3^3$ and
thus by necessity (since the Bockstein spectral sequence must
eventually converge to zero), $B_2^4=B_3^4$. Thus in particular,
$B_3^* \neq 0$ in positive dimensions. Thus we conclude:
$\expo(\bar{H}^*(G(\mathfrak{sl}_2) ;\mathbb{Z})) > p^2$. On the other
hand, $G(\mathfrak{sl}_2)$ contains $G(\mathfrak{S}) \times
\mathbb{Z}/p\mathbb{Z}$ as a subgroup of index $p$. Thus by transfer
arguments, 
$$
\expo(\bar{H}^*(G(\mathfrak{sl}_2) ;\mathbb{Z})) \leq 
p\expo(\bar{H}^*(G(\mathfrak{S});\mathbb{Z})) = p^3.
$$
Thus
$$
\expo(\bar{H}^*(G(\mathfrak{sl}_2) ;\mathbb{Z})) = p^3.
$$
Note this is neither the order nor the exponent of
$G(\mathfrak{sl}_2)$
as $G(\mathfrak{sl}_2)$ has order $p^6$ and exponent $p^2$.

The authors would like to thank T. Weigel for many useful discussions.

\end{document}